\documentclass[11pt]{amsart}
\usepackage{hyperref}
\hypersetup{nesting=true,debug=true,naturalnames=true}
\usepackage{graphicx,amssymb,upref}

\usepackage{amsmath,amsthm}
\usepackage{amsfonts}
\usepackage{latexsym}
\usepackage{amsbsy}
\usepackage{amssymb,mathscinet}

\newtheorem{theorem}{Theorem}[section]
\newtheorem{corollary}[theorem]{Corollary}
\newtheorem{lemma}[theorem]{Lemma}
\newtheorem{proposition}[theorem]{Proposition} 

\theoremstyle{definition}
\newtheorem{definition}[theorem]{Definition}
\newtheorem{remark}[theorem]{Remark}
\newtheorem{notation}[theorem]{Notation}
\newtheorem{Example}[theorem]{Example}
\theoremstyle{remark}

\numberwithin{equation}{section}

\newcommand{\be}{\mathbb E}
\newcommand{\bn}{\mathbb N}

\newcommand{\ot}{\otimes}
\newcommand {\id} {{\textrm{id}}}
\newcommand{\wt}{\widetilde}
\newcommand{\wT}{\wt{T}}

\newcommand{\wA}{\wt{A}}

\newcommand{\wS}{\wt{S}}

\begin{document}

	\title[Wold-type decomposition for doubly twist. left-inv. cov. repr.]
	{Wold-type decomposition for doubly twisted left-invertible covariant representations}

	\date{\today}

	
		\author[Kumar]{Niraj Kumar\textsuperscript{*}}
	\address{The LNM Institute of Information Technology, Rupa ki Nangal, Post-Sumel, Via-Jamdoli
		Jaipur-302031,
		(Rajasthan) INDIA}
	\email{24pmt002@lnmiit.ac.in}

\author[Rohilla]{Azad Rohilla}
	
	\address{Sanskaram University, Patauda, Jhajjar, Haryana, 124108,
		 INDIA}
	\email{18pmt005@lnmiit.ac.in, azadrohilla23@gmail.com}
	
	\author[Trivedi]{Harsh Trivedi}
	
	\address{The LNM Institute of Information Technology, Rupa ki Nangal, Post-Sumel, Via-Jamdoli
		Jaipur-302031,
		(Rajasthan) INDIA}
	\email{harsh.trivedi@lnmiit.ac.in, trivediharsh26@gmail.com}

	\thanks{*corresponding author}

	
	\begin{abstract}
		In this article, we have introduced the notion of a near-isometric covariant representation of a $C^*$-correspondence. The other objective is to provide a unified approach to several known results for a large class of left-invertible covariant representations of a product system and prove Wold-type decomposition for the case of doubly twisted left-invertible covariant representations and study some applications. 
	\end{abstract}


	\maketitle
	\section{Introduction}
	
 Wold \cite{W} showed that an isometry on a Hilbert space can be decomposed orthogonally as a shift and a unitary uniquely as a direct summand. 
S{\l}oci{\'n}ski \cite{Sl80} provided this decomposition for doubly commuting case. This decomposition and its applications to weakly stationary random fields is also discussed in a book by Mandrekar and Redett \cite{VR17}. 
In recent years, there has been several results in this direction, see \cite{AS25, MS99,Po89,S01, S14,SZ08,Sl80,TV19}. In \cite{La22}, Lata, Pokhriyal and Singh proved the orthogonal decomposition for near-isometry and also explored the case of doubly commuting. In \cite{Ra22} Rakshit, Sarkar and Suryawanshi extended the notion of doubly commuting isometries and introduced the notion of $\mathcal{U}_n$-twisted isometries and proved that they admit the Wold-type decomposition.

	 Wold-type decomposition for Cuntz (row) isometries \cite{C77} was studied by Popescu \cite{Po89}. Pimsner \cite{P97} using covariant representations generalized the notion of algebras generated by Cuntz isometries and proved that the covariant representation of a $C^*$-correspondence is isometric is equivalent to that the representation extends to the Toeplitz algebra of the correspondence. Muhly and Solel in \cite{MS99} proved the Wold-type decomposition by considering isometric covariant representations. 
	Doubly comuting case of this result is due to Skalski-Zacharias \cite{SZ08}.

    Popescu \cite{PG23} studied \(k\)-tuples of doubly twisted (\(U\)-commuting) row isometries, where \(U\) is a family of commuting unitary operators, and established a Wold-type decomposition. Recently, in \cite{SS26} Wold-type decomposition and as an application unitary extension for twisted isometric covariant representation associated with product system were studied by Solel and Suryawanshi. The present paper considers covariant representations on product systems that are left-invertible, and which satisfy doubly twisted commutation relations. 

    In this paper, in Section \ref{secnr}, we obtained Wold-type decomposition for a near-isometric covariant
representation of a $C^*$-correspondence. In Section \ref{Secdc}, we studied doubly twisted left-invertible covariant representations and obtained some basic results. In Section \ref{secwd}, Wold-type decomposition for several doubly twisted left-invertible representations of a product system is proved and our approach for the decomposition is based on a unified approach for doubly commuting left-invertible operators given recently in \cite{NZ25}.
	

	\subsection{Preliminaries and Notations}
	 In this subsection, we recall some definitions and results related to Hilbert $C^*$-modules and covariant representations of
	$C^*$-correspondences (cf. \cite{La95, P97,MR1648483}) and recall the Wold-type decomposition due to Muhly and Solel \cite{MS99} for isometric covariant representations.
    
	Let $\mathcal A$ be a $C^*$-algebra, and $F$ a Hilbert $\mathcal A$-module. The set of all adjointable operators on $F$ is denoted by $\mathcal L(F).$
	 The module $F$ is said to be a {\it
		$C^*$-correspondence over $\mathcal A$} if $F$ has a left action  by a non-zero $*$-homomorphism
	$\pi:\mathcal A\to \mathcal L(F)$ such that
	\[
	a\eta:=\pi(a)\eta \quad \quad (a\in\mathcal A, \eta\in F).
	\]
	The map $\pi$ is {\it essential} in the sense that 
	$F=\bigvee\pi(\mathcal A)F$ and will be used throughout. A $C^*$-correspondence inherits the operator space structure induced by considering it as a corner in the relevant linking algebra. For a $C^*$-correspondence $E$ over $\mathcal A,$ we can construct the balanced tensor product $E\otimes_{\pi} F$ satisfying
	\[
	(\gamma_1 a)\otimes \zeta_1=\gamma_1\otimes \pi(a)\zeta_1,
	\]
	\[
	\langle\gamma_1\otimes\zeta_1,\gamma_2\otimes\zeta_2\rangle=\langle\zeta_1,\pi(\langle\gamma_1,\gamma_2\rangle)\zeta_2\rangle,
	\]
	for every $\zeta_1,\zeta_2\in F;$ $\gamma_1,\gamma_2\in E$ and $a\in\mathcal
	A.$  We simply denote $E\otimes_{\pi} F$ by  $E\otimes F.$  
	\begin{definition}
		Consider a Hilbert space $\mathcal{K},$ and $E$ a $C^*$-correspondence  over $\mathcal A.$ Then for a representation
		$\sigma:\mathcal A\to B(\mathcal{K})$ and a linear map $A:
		E\to B(\mathcal{K}),$ an ordered pair $(\sigma, A)$ is
		said to be a {\it covariant representation} of $E$ on $\mathcal{K}$ (cf. \cite{MR1648483}) if
		\[
		\sigma(a)A(\eta)\sigma(b)=A(\pi(a)\eta b) \quad \quad (\eta\in E;
		a,b\in\mathcal A).
		\]
		If $A$ is completely bounded (respectively, contractive), then we call $(\sigma, A)$ a {\it completely bounded} (respectively,  {\it
			contractive}). Furthermore, it is said to be {\it isometric} or {\it Toeplitz} if
		\[
		\sigma(\langle \eta,\zeta\rangle)=A(\eta)^*A(\zeta) \quad \quad
		(\eta,\zeta\in E).
		\]
	\end{definition}
    Define $\widetilde{A}:~\mbox{$E\otimes_{\sigma} \mathcal{K}\to \mathcal{K}
			$}$ by
		\[
		\widetilde{A}(\eta\otimes h):=A(\eta)h \quad \quad (\eta\in E,
		h\in\mathcal{K}).
		\]
		Then $\widetilde{A}(\pi(a)\otimes I_{\mathcal{K}
			})=\sigma(a)\widetilde{A}$, $a\in\mathcal A$.
	
	\begin{lemma}{\rm (\cite[Lemma
			3.5]{MR1648483})}\label{MSL}
		The map $(\sigma, A)\mapsto \widetilde A$ gives a bijection
		between the set of each completely bounded (respectively, contractive), covariant
		representations $(\sigma, A)$ of $E$ on $\mathcal{K}$ and the
		set of each bounded (respectively, contractive) linear maps
		$\wA$ satisfying $\widetilde{A}(\pi(a)\otimes I_{\mathcal
			K})=\sigma(a)\widetilde{A}$, $a\in\mathcal A$. Indeed, the condition that $(\sigma, A)$ is isometric is equivalent to $\widetilde
		A$ being isometry.
	\end{lemma}
	If $\widetilde{A}$ is co-isometry, then the covariant representation $(\sigma, A)$ is said to be {\it fully co-isometric} . Consider a $C^*$-correspondence $E$ over $\mathcal A.$ Then for every $n\in \mathbb{N}$,
	$E^{\otimes n}: =E\otimes_{\pi} \cdots \otimes_{\pi}E$ ($n$ fold tensor product) is a $C^*$-correspondence and the left action is given by  $$\pi_n(a)(\eta_1 \otimes \cdots \otimes \eta_n):=\pi(a)\eta_1\otimes \cdots \otimes\eta_n.$$   We use $E^{\otimes 0}:=\mathcal{A}$ and $\mathbb Z_+ :=\mathbb N\cup\{0\}$. 
	For each $ n \in \mathbb{N},$ and for a given covariant representation $(\sigma, A),$ define the map $\wA_n: E^{\ot n}\otimes_{\sigma} \mathcal{K}\to \mathcal{K}$ by
	\[ \wA_n (\eta_1 \ot \cdots \ot \eta_n \ot h) = A (\eta_1) \cdots A(\eta_n) h,\]
	where $\eta_1, \ldots, \eta_n\in E, h \in \mathcal{K}$. Therefore
	$\wA_n=\wA(I_{E} \otimes \wA) \cdots (I_{E^{\otimes n-1}} \otimes  \wA ).$

	The {\it Fock  $C^*$-correspondence} $\mathcal{F}(E):= \bigoplus_{n \in \mathbb Z_+}E^{\otimes n}$ where the left action  of $\mathcal{A}$ on $\mathcal{F}(E)$ is denoted by $\pi_{\infty}: \mathcal{A} \longrightarrow L(\mathcal{F}(E))$ and given as $$\pi_{\infty}(a)(\oplus_{n \in \mathbb Z_+}\eta_n):=\oplus_{n \in \mathbb Z_+}\pi_{n}(a)\eta_n , \:\: \eta_n \in E^{\otimes n}.$$
	For $\eta \in E$, the {\it creation operator} $A_{\eta}$ on $\mathcal{F}(E)$ is defined by $$A_{\eta}(\tau):=\eta \otimes \tau, \:\: \tau \in E^{\otimes n}, n \in \mathbb Z_+.$$

	\begin{definition}
		Suppose that $E$ be a $C^*$-correspondence over $\mathcal A,$ and let $\psi $ be a $*$-homomorphism of $\mathcal{A}$ on $\mathcal{K}$  Then the {\it induced representation} (cf. \cite{R74}) is an isometric representation $(\rho, S)$ of $E$ on $\mathcal{F}(E)\otimes_{\psi}\mathcal{K}$ and is defined by
		\begin{align*}
			\rho(a):&=\pi_{\infty}(a) \otimes I_{\mathcal{K}} \:\:, a \in \mathcal{A}&\\
			S(\eta):&=A_{\eta}\otimes I_{\mathcal{K}} ,\:\: \eta \in E.
		\end{align*}
	\end{definition}

	\begin{definition}	
		\begin{itemize}
			\item[(i)] For a completely bounded, covariant representation $(\sigma, A)$ of $E$ on $\mathcal{K},$ a closed subspace  $\mathcal S$ of $ \mathcal{K}$ is $(\sigma, A)$-{\it invariant} (respectively,  $(\sigma, A)$-{\it reducing}) (cf. \cite{SZ08}) if it   is $\sigma(\mathcal A)$-invariant, and if $\mathcal S$ (respectively,  both $ \mathcal{S}, \mathcal{S}^{\bot}$) is invariant under each operator $A (\eta)$ for $\eta \in E.$ Then the canonical restriction provides a new representation of $E$ on $\mathcal{S}$ and it will denoted by $(\sigma, A)|_{\mathcal{S}}.$
			
			\item[(ii)]  For a closed $\sigma(\mathcal{A})$-invariant subspace $\mathcal{S},$ we define $$\mathfrak{A}_{n}(\mathcal{S}):=\bigvee \{\wA_n(\eta_1\ot\eta_2\ot\cdots \ot\eta_n\ot h) \: : \: \eta_i\in E , h \in \mathcal{S} \},$$ for $n\in\mathbb N$ and $\mathfrak{A}_{0}(\mathcal{S}):=\mathcal{S}$. If $\mathcal{S}\perp\mathfrak{A}_{n}(\mathcal{S}),$ for each $n \in \mathbb{N},$ then $\mathcal{S}$ is said to {\it wandering}. Moreover, $(\sigma, A)$ has {\it generating wandering subspace property} if
			$$\mathcal{K}=\bigvee_{n\geq0}\mathfrak{A}_n(\mathcal{S}).$$
		\end{itemize}

	\end{definition}
	
	
	\begin{theorem}\cite{MS99}\label{ThmMS}
		Every isometric covariant representation $(\sigma, A)$ of $E$ on $\mathcal{K}$ decomposes into a direct sum $(\sigma_1, A_1)\bigoplus(\sigma_2, A_2)$ on $\mathcal{K}=\mathcal{K}_1 \bigoplus \mathcal{K}_2$ such that $(\sigma_1, A_1)=(\sigma, A)|_{\mathcal{K}_1}$ is an induced representation and $(\sigma_2, A_2)=(\sigma, A)|_{\mathcal{K}_2}$ is fully coisometric. This decomposition is unique, that is, if $\mathcal{S}$ reduces $(\sigma, A),$  and if the restriction $(\sigma, A)|_{\mathcal{S}}$ is induced (respectively,  fully coisometric), then $\mathcal{S} \subseteq \mathcal{K}_1$(respectively,  $\mathcal{S} \subseteq \mathcal{K}_2$).  Moreover, $\mathcal{K}_1:=\bigoplus_{n \in \mathbb Z_+}\mathfrak{A}_{n}(\mathcal{N}(\wA^*)),$ and hence	
		$$\mathcal{K}_2:=\left(\bigoplus_{n \in \mathbb Z_+}\mathfrak{A}_{n}(\mathcal{N}(\wA^*))\right)^{\bot}=\bigcap_{n\in  \mathbb Z_+}\text{ran} (\widetilde{A}_n),$$ where $A_0:=\sigma.$ 
	\end{theorem}
	
	The previous result was extended in the following results (see \cite{HS19, RTV23}) which is a generalization of Wold-type decomposition for concave operators (see \cite{R88,S01,O05}):

\begin{theorem}\label{MT1}
Suppose $(\sigma, A)$ is a completely bounded, left-invertible covariant representation of $E$ on  $\mathcal{K}$ such that any one of the following conditions hold:
\begin{enumerate}
\item[(1)]  $(\sigma, A)$ is concave, that is, \\ $\| \wt{A}_2(\eta \otimes h)\|^2+\|\eta \otimes h\|^2 \leq 2 \|(I_E \otimes  \wt{A})(\eta \otimes h)\|^2, \eta \in E^{\otimes 2}, h \in \mathcal{K}.$ \cite{HS19}
\item[(2)]  $\|(I_{E} \otimes \wt{A})(\zeta)+\kappa\|^2 \leq 2(\|\zeta\|^2+\|\wt{A}(\kappa)\|^2), \:\:\:\: \zeta \in E^{\otimes 2} \otimes \mathcal{K},\kappa \in E \otimes \mathcal{K}.$ \cite{HS19}
\item [(3)] $(\sigma, A)$ is completely contractive satisfying  $\|\wt{A}(\xi)\|^2+\|\wt{A}^*_2\wt{A}(\xi)\|^2 \leq 2\|\wt{A}^*\wt{A}(\xi )\|^2, \:  \xi \in E \otimes \mathcal{K}.$ \cite{HS19}
\item [(4)] Suppose $ l>0 $ and $(l_m) $ is a sequence of non-negative numbers such that $$\|\widetilde{A}_m (\xi_{m} )\|^2 \leq l_m \big(\| (I_{E^{\ot {m-1}}} \ot \widetilde{A}) (\xi_{m})\|^2 -\| (\xi_{m})\|^2\big) + l\|(\xi_{m})\|^2$$
where $  \xi_{m} \in E^{\ot m} \ot \mathcal{K}, m\in \mathbb{N}.$\cite{RTV23}

\end{enumerate}
Then $(\sigma, A)$ admits  {\it Wold-type decomposition} (when $l=1$ in condition $(4)$). This decomposition is unique and both the reducing subspaces are given by
$$\mathcal{K}_1=\bigvee_{m\geq 0}\mathfrak{A}_m(\mathcal{N}(\wA^*)), \:\:\:\:\: \mathcal{K}_2=\bigcap_{m \geq 1}\wt{A}_m(E^{\otimes m} \otimes \mathcal{K}).$$   Indeed, if $(\sigma, A)$ is analytic (pure), that is, $\mathcal{K}_{2}=\{0\}$, then $(\mathcal{N}(\wA^*))$ is the generating wandering subspace for $(\sigma, A),$ that is, $\mathcal{K}=\bigvee_{m\geq 0}\mathfrak{A}_m(\mathcal{N}(\wA^*)).$
\end{theorem}

In \cite{NZ25}, the Wold-type decomposition from \cite{O05,R88,S01} is studied in the doubly commuting setting. In this paper we will extend doubly commuting Wold-type decomposition \cite{NZ25} for the doubly twisted covariant representation case. 

\section{Wold-type decomposition for near-isometric covariant representations of C*-correspondences}\label{secnr}
	
	In this section, we introduce the notion of a near-isometric covariant representation and prove a Wold-type decomposition based on \cite{La22}.
	
{\bf{Nagy-Foias-Langer decomposition}:} Every contraction $T$ on  $\mathcal{K}$ (\cite{NF70}) decomposes $\mathcal{K} = \mathcal{K}_0\oplus \mathcal{K}_1,$ where $\mathcal{K}_0$ and $\mathcal{K}_1$ be two subspaces of $\mathcal{K}$ reduce $T$ such that the restriction of $T$ to $\mathcal{K}_0$ is unitary and the restriction to $\mathcal{K}_1$ is completely non unitary.

\begin{definition}
    A completely contractive covariant representation $(\sigma,A)$ is said to have a Wold-type decomposition if $\mathcal{K}$ admits an orthogonal decomposition $\mathcal{K} = \mathcal{K}_0\oplus \mathcal{K}_1,$ where $\mathcal{K}_0$ and $\mathcal{K}_1$ are two $(\sigma,A)$-reducing subspaces of $\mathcal{K}$ and $(\sigma,A)|_{\mathcal{K}_0}$ is an isometric as well as fully co-isometric representation and $(\sigma,A)|_{\mathcal{K}_1}$ is completely non unitary.
\end{definition}
 
\begin{theorem}(see \cite[Theorem 1.5]{AD26}).
  If $(\sigma , A)$ is a completely contractive covariant representation of $E$ on  $\mathcal{K}.$ Then $\mathcal{K}$ admits a unique orthogonal decomposition $\mathcal{K} = \mathcal{K}_0\oplus \mathcal{K}_1$ we have the following:
  \begin{enumerate}
      \item $\mathcal{K}_0$ and $\mathcal{K}_1$ reduce $(\sigma , A).$
      \item The restriction $(\sigma,A)|_{\mathcal{K}_0}$ is an isometric as well as fully co-isometric representation.
      \item The restriction $(\sigma , A)|_{\mathcal{K}_1}$ is completely non unitary.
  \end{enumerate}
\end{theorem}
	
	\begin{definition}	
		Let $(\sigma, A)$ be a completely bounded covariant representation of $E$ on  $\mathcal{K}.$ The representation $(\sigma, A)$ is said to be {\it near-isometric} if the following two conditions hold:
		\begin{itemize}
			\item[(a)] \label{a} For each $ \zeta\in E \otimes \mathcal{K},$ we have  $\alpha \|\zeta\| \le\|\wt{A}\zeta\| \le \|\zeta\|,$ for some $\alpha> 0$.
			\item[(b)] For every $m\ge 1$ we have  $\text{ran}(\wt{A}^*_m\wt{A}_{m+1})$ $\subseteq {E^{\otimes m} \otimes \text{ran}\wt{A}}$.

		\end{itemize}
		
	\end{definition}

\begin{Example}
     Suppose $(\rho, V)$ is an isometric covariant representation and fix a constant $0<\beta\le1.$ Define a covariant representation $(\sigma, A)$ of $E$ on $\mathcal{K}$ defined by $\sigma(a):=\rho(a)$ and $A(x):=\beta V(x)$ for each $a\in\mathcal A$ and $x\in E.$  Then for each $\zeta\in E\otimes\mathcal{K}$,
$\|\widetilde A\zeta\|=\beta\|\zeta\|.$ Thus,
$\beta\|\zeta\|\le\|\widetilde A\zeta\|\le\|\zeta\|.$ For every $m\ge1$,
$\widetilde A_m=\beta^m \widetilde V_m$ and we get\begin{align*}
\widetilde A_m^* \widetilde A_{m+1}
= \beta^{2m} (I_{E^{\otimes m}}\ot\widetilde A).
\end{align*}
 It follows that 
\begin{equation*}
   \text{ran}(\widetilde A_m^*\widetilde A_{m+1})
\subseteq
E^{\otimes m}\otimes \operatorname{ran}\widetilde A,  
\end{equation*} and hence $(\sigma,A)$ is a near-isometric covariant representation of $E$.
\end{Example}

\begin{Example}
     In this example we shall use several properties of bilateral shift discussed in \cite[Subsection 7.2]{RTV23}. Consider $l^2 (\mathbb{Z})$ with the canonical orthonormal basis $\{e_m\}_{m\in \mathbb{Z}},$ and let $J_n = \{1,\ldots,n\}$. For each $i\in J_n,$ let $\{w_{i,m}\}_{m\in \mathbb{Z}}\subset \mathbb{C}$ be a bounded weight sequence. Define operators $V_i\in B(l^2(\mathbb{Z}))$ by $V_i(e_m)=w_{i,m}e_{i+nm},m\in \mathbb{Z}.$ Let $E=\mathbb{C}^n$ with orthonormal basis $\{\delta_i\}_{i\in J_n}.$ Define a covariant representation $(\rho,S^w)$ of $E$ on $l^2(\mathbb{Z})$ by $\rho(\lambda)=\lambda I_{l^2(\mathbb{Z})},S^w(\delta_i)=V_i,\lambda\in\mathbb{C},i\in J_n.$
     For $\lambda\in\mathbb{C}$ and  we have $S^w(\lambda\delta_i)=\lambda S^w(\delta_i)=\lambda V_i=\rho(\lambda)V_i=\rho(\lambda)S^w(\delta_i),$ and similarly $S^w(\delta_i \lambda)=S^w(\delta_i)\rho(\lambda).$ Hence, $(\rho,S^w)$ is a covariant representation of $E$ on $l^2(\mathbb{Z}).$\\
    Assumptions on the weights: $\exists \gamma>0$ such that $\gamma\le |w_{i,m}|\le 1,\forall i\in J_n$ and all $m\in \mathbb{Z}\setminus\{0\},$ and $w_{i,0}=0,\forall i\in J_n.$ Now we show that  $S^w$ is a near-isometry. For that we verify the following two conditions : (a) let $\wS^{w}:E\ot l^2(\mathbb{Z})\to l^2(\mathbb{Z})$ be the operator associated with $S^w.$ For all $\delta_i\ot e_m\in E\ot l^2(\mathbb{Z}),$ we have $\wS^{w}(\delta_i\ot e_m)=V_ie_m,$ using the orthogonality of the ranges of $\{V_i\}_{i=1}^n$ and the bounds on the weights, we get $\gamma\|x\|\le \|\wS^{w}x\|\le \|x\|$ for all $x\in E\ot l^2(\mathbb{Z}).$ \\
    (b) Note that 
     \begin{align*}
        \langle \wS^{w^*}\wS^{w}(\delta_{i}\ot e_m),(\delta_{q}\ot e_k)\rangle & = \langle V_i e_m,V_qe_k\rangle  \\&
        = w_{i,m}\overline{w_{q,k}} \langle e_{i+nm},e_{q+nk}\rangle \\&
        =  w_{i,m}\overline{w_{q,k}} \langle \delta_{i},\delta_{q}\rangle \langle e_{m},e_{k}\rangle\\&
        = w_{i,m}\overline{w_{q,k}}\langle \delta_{i}\ot e_m ,\delta_{q}\ot e_k\rangle.
    \end{align*}
     For each $l\ge 1,$ we have
    \begin{align*}
        \wS_l^{w}(\delta_{i_l}\ot\ldots \ot\delta_{i_1}\ot e_m) & = S^w(\delta_{i_l})\ldots S^w(\delta_{i_1})e_m\\&
        = w_{i_1,m}\left(\Pi_{j=2}^l w_{i_j,(\sum_{k=0}^{j-2}i_{l-k-1}n^k+n^{j-1}m)} \right)e_{(\sum_{k=0}^{l-1}i_{l-k}n^k+n^{l}m)}.
    \end{align*}
    Therefore

    \begin{align*}
        & \langle\wS_l^{{w}^*}\wS_l^{w}(\delta_{i_l}\ot\ldots \ot\delta_{i_1}\ot e_m),(\delta_{r_l}\ot\ldots \ot\delta_{r_1}\ot e_t)\rangle \\& 
        = \langle w_{i_1,m}\left(\Pi_{j=2}^l w_{i_j,(\sum_{k=0}^{j-2}i_{l-k-1}n^k+n^{j-1}m)} \right)e_{(\sum_{k=0}^{l-1}i_{l-k}n^k+n^{l}m)}, \\&      \quad  \quad  \quad  \quad w_{r_1,t}\left(\Pi_{j=2}^l w_{r_j,(\sum_{k=0}^{j-2}r_{l-k-1}n^k+n^{j-1}t)} \right)e_{(\sum_{k=0}^{l-1}r_{l-k}n^k+n^{l}t)}\rangle\\& 
        = w_{i_1,m}\left(\Pi_{j=2}^l w_{i_j,(\sum_{k=0}^{j-2}i_{l-k-1}n^k+n^{j-1}m)} \right)\overline{w_{r_1,t}\left(\Pi_{j=2}^l w_{r_j,(\sum_{k=0}^{j-2}r_{l-k-1}n^k+n^{j-1}t)} \right)}\\& 
        \quad  \quad  \quad  \quad \langle e_{(\sum_{k=0}^{l-1}i_{l-k}n^k+n^{l}m)},e_{(\sum_{k=0}^{l-1}r_{l-k}n^k+n^{l}t)}\rangle\\& 
        = w_{i_1,m}\left(\Pi_{j=2}^l w_{i_j,(\sum_{k=0}^{j-2}i_{l-k-1}n^k+n^{j-1}m)} \right)\overline{w_{r_1,t}\left(\Pi_{j=2}^l w_{r_j,(\sum_{k=0}^{j-2}r_{l-k-1}n^k+n^{j-1}t)} \right)} \\&
        \quad  \quad  \quad  \quad \langle \delta_{i_1},\delta_{r_1}\rangle, \ldots,\langle \delta_{i_l},\delta_{r_l}\rangle\langle e_m,e_t\rangle\\&
        = w_{i_1,m}\left(\Pi_{j=2}^l w_{i_j,(\sum_{k=0}^{j-2}i_{l-k-1}n^k+n^{j-1}m)} \right)\overline{w_{r_1,t}\left(\Pi_{j=2}^l w_{r_j,(\sum_{k=0}^{j-2}r_{l-k-1}n^k+n^{j-1}t)} \right)}\\&
        \quad  \quad  \quad  \quad \langle (\delta_{i_l}\ot\ldots \ot\delta_{i_1}\ot e_m),(\delta_{r_l}\ot\ldots \ot\delta_{r_1}\ot e_t)\rangle.
    \end{align*}
      Consequently, $ran(\wS_m^{w^*} \wS_{m+1}^{w} )\subseteq E^{\ot m}\ot \text{ran}\wS_m^{w} .$
     Under the above assumptions  on the weights, the bilateral shift representation $(\rho,S^w)$ is a covariant representation of $E$ that satisfies both conditions of near-isometric. 
      \end{Example}

	\begin{theorem}\label{ThmW}
		Let $(\sigma , A)$ be a near-isometric covariant representation of $E$  on  $\mathcal{K}$. Then we can decompose the representation $(\sigma, A)$ as $(\sigma_1, A_1)\bigoplus(\sigma_2, A_2)$ on $\mathcal{K}=\mathcal{K}_1 \bigoplus \mathcal{K}_2$ such that $(\sigma_1, A_1):=(\sigma, A)|_{\mathcal{K}_1}$ is an induced representation and $(\sigma_2, A_2):=(\sigma, A)|_{\mathcal{K}_2}$ an invertible covariant representation. This decomposition is uniquely determined if $S$  reduces $(\sigma, A),$  and if the restriction $(\sigma, A)|_{S}$ is induced (respectively, invertible), then $S \subseteq \mathcal{K}_1$(respectively,  $S \subseteq \mathcal{K}_2$).  Moreover, $\mathcal{K}_1:=\bigoplus_{n \in \mathbb Z_+}\mathfrak{A}_{n}(\mathcal{N}(\wt{A}^*)),$ and hence	
		$$\mathcal{K}_2:=\bigcap_{n\in  \mathbb Z_+}\text{ran}(\wt{A}_n),$$ where $A_0:=\sigma.$
	\end{theorem}
	\begin{proof}: We know that  $\mathcal{N}(\wt{A}^*)=\mathcal{K}\ominus\widetilde{A}(E\otimes \mathcal{K})$. We are given that $(\sigma, A)$ is a near-isometric covariant representation, therefore $\text{ran}(\wt{A}^*_m\wt{A}_{m+1})$ $\subseteq {E^{\otimes m} \otimes ran(\wt{A}})$ for each $m\ge 0$. Hence, for each $\xi_1,\xi_2,\ldots,\xi_m$ and $\eta_1,\eta_2,\ldots,\eta_{m+1}\in E $, $w\in \mathcal{N}(\wt{A}^*), h\in \mathcal{K}$ we have 
		\begin{align*}
			& \langle A(\xi_1),A(\xi_2) \cdots A(\xi_m)w, A(\eta_1)A(\eta_2) \cdots A(\eta_{m+1})h\rangle
			\\& =\langle \xi_1\ot\xi_2 \ot \cdots \ot\xi_m\ot w, \wt{A}^*_m\wt{A}_{m+1}(\eta_1\ot\eta_2 \ot \cdots \ot\eta_{m+1}\ot h)\rangle
			\\& =0.
		\end{align*}
		So, $\mathfrak{A}_{m}(\mathcal{N}(\wt{A}^*))\perp \text{ran}(\wt{A}_{m+1})$ for each $m\ge0$ and therefore $\mathfrak{A}_{m}(\mathcal{N}(\wt{A}^*))\perp \mathfrak{A}_{k}(\mathcal{N}(\wt{A}^*)), m\neq k$. Now using $\widetilde{A}$ is bounded below, it implies that $\mathfrak{A}_{m}(\mathcal{N}(\wt{A}^*))$ is a closed subspace of $\mathcal{K}$.  Set $\mathcal{K}_1:=\bigoplus_{n \in \mathbb Z_+}\mathfrak{A}_{n}(\mathcal{N}(\wt{A}^*)),$ and therefore $$\mathcal{K}_2:=\bigcap_{n\in  \mathbb Z_+}\text{ran}(\wt{A}_n).$$
		We want to show that $\mathcal{K}_1^{\perp}=\mathcal{K}_2$.\
		Since  $\mathfrak{A}_{m}(\mathcal{N}(\wt{A}^*))\perp \text{ran}(\wt{A}_{m+1})$ for each $m\ge0$, we have $\mathcal{K}_1\perp \mathcal{K}_2$ which implies that $\mathcal{K}_2\subseteq \mathcal{K}_1^\perp$. For the other inclusion, note that $\mathfrak{A}_{m}(\mathcal{N}(\wt{A}^*))\perp \text{ran}(\wt{A}_{m+1})$ for each $m\ge0$, which gives us $\text{ran}(\wt{A}_m)=\mathfrak{A}_{m}(\mathcal{N}(\wt{A}^*))\bigoplus \text{ran}(\wt{A}_{m+1})$ for each $m\ge0$. Therefore, if $x\perp \mathfrak{A}_{m}(\mathcal{N}(\wt{A}^*))$ for each $m\ge0$, then $x\in \text{ran}(\wt{A}_m)$, for each $m\ge0$. Hence $\mathcal{K}_2=\mathcal{K}_1^{\perp}$.

		\end{proof}
        
        \section{Basic notations and preparatory results: doubly twisted left-invertible covariant representations}\label{Secdc}
	The tensor product system of Hilbert spaces was studied by Arveson \cite{A89} and the discrete product system of $C^*$-correspondences was studied by Fowler \cite{F02}. In this section, we recapture the several notions related to the product systems of $C^*$-correspondences and the twisted doubly commuting case from (see \cite{F02,S06,  SZ08, SS26}) and examine wandering subspaces for the doubly twisted covariant representations of product systems. 
	
	\begin{definition}
		\begin{itemize}
			\item[(1)]  Let $q\in \mathbb N.$ A {\it product system} (cf. \cite{F02})
			$\be$ on $\mathbb Z^q_+$ is defined as a family of $C^*$-correspondences $\{E_1, \ldots,E_q\},$ along with the unitary
			isomorphisms $u_{i,k}: E_i \ot E_k \to E_k \ot E_i$ ($i>k$). Using these isomorphisms, for each
			${\bf m}=(m_1, \cdots, m_q) \in \mathbb Z^q_+$ the correspondence $\be ({\bf m})$ is identified with $E_1^{\ot^{ m_1}} \ot \cdots \ot E_q^{\ot^{m_q}}.$ We define maps $u_{i,i} := \id_{E_i \ot E_i}$ and $u_{i,k} := u_{k,i}^{-1}$ when $i<k.$
			
			\item[(2)]  When $i<k,$ let $\{U_{ik}\}$ be $\begin{pmatrix} q \\ 2 \end{pmatrix}$ commuting unitaries on $\mathcal{K}$ such that $U_{ki}:=U^*_{ik}.$ Then call $\{U_{ik}\}_{i<k}$ a {\it twist}. Assume $\be$ to be a product system over $\mathbb Z^q_+$. A {\it completely bounded, twisted covariant representation}  of $\be$ on $\mathcal{K}$ is defined as a tuple $(\sigma, A^{(1)}, \ldots, A^{(q)})$, where $\sigma$ is a
			representation of the $C^*$-algebra $\mathcal A$ on $\mathcal{K}$,  $A^{(i)}:E_i \to B(\mathcal{K})$ are  completely bounded linear  maps satisfying
			\[ A^{(i)}(a \xi_i b) = \sigma(a) A^{(i)}(\xi_i) \sigma(b), \;\;\; a,b \in \mathcal A, \xi_i \in E_i,\]
			and satisfy the twisted commutation relations (see \cite{SS26})
			\begin{equation} \label{twist} \wA^{(i)} (I_{E_i} \ot \wA^{(j)}) = U_{ij}\wA^{(j)} (I_{E_j} \ot \wA^{(i)})(u_{i,j}\ot I_\mathcal{K}) \end{equation} and
            \begin{equation} \label{ucomm} \wA^{(l)}(I_{E_l} \otimes U_{ij})=U_{ij} \wA^{(l)}, U_{ij}\in \sigma(\mathcal{A})',\end{equation}
			with $1\leq i,j,l\leq q$ and $i\neq j$.
            Moreover, the completely bounded covariant representation $(\sigma, A^{(1)}, \ldots, A^{(q)})$ is called twisted {\it isometric} if each $(\sigma, A^{(i)})$ is isometric as a covariant representation of the $C^*$-correspondence $E_i$, and similarly we can define the notion of twisted {\it near-isometric}.
		\end{itemize}
	\end{definition}
    
	\begin{remark}
    Equation (\ref{twist}) admits an equivalent representation given by 
             \begin{align}\label{eqn a}
              \wA^{(i)} (I_{E_i} \ot \wA^{(j)}) &= U_{ij}\wA^{(j)} (I_{E_j} \ot \wA^{(i)})(u_{i,j}\ot I_\mathcal{K}) \\& 
            = \wA^{(j)}(I_{E_j}\ot  U_{ij}) (I_{E_j} \ot \wA^{(i)}) (u_{i,j}\ot I_\mathcal{K}) \nonumber\\&
            = \wA^{(j)}(I_{E_j}\ot  U_{ij}\wA^{(i)})(u_{i,j}\ot I_\mathcal{K})\nonumber\\& 
            = \wA^{(j)}(I_{E_j}\ot  \wA^{(i)}(I_{E_i}\ot  U_{ij}))(u_{i,j}\ot I_\mathcal{K})\nonumber\\& 
            = \wA^{(j)}(I_{E_j}\ot  \wA^{(i)})(I_{E_j\ot E_i}\ot  U_{ij})(u_{i,j}\ot I_\mathcal{K})\nonumber\\& 
            = \wA^{(j)}(I_{E_j}\ot  \wA^{(i)})(u_{i,j}\ot U_{ij})\nonumber.
      \end{align}
    \end{remark}

\begin{notation}
    We denote $\mathbb{E}$ for the product system over $\mathbb Z^k_+,k\in \mathbb{N}$ in the following upcoming result.
\end{notation}

    \begin{definition}  \label{dcom}
		A completely bounded, covariant representation $(\sigma, A^{(1)}, \ldots, A^{(k)})$ of  $\be$ on $\mathcal{K}$ is said to be {\it doubly twisted}  with respect to the twist $U_{ij}$ if it is twisted and for each $i\neq j \in J_k,$ we have \begin{equation}\label{doubly}\wA^{(j)^*} \wA^{(i)} =  (I_{E_j} \ot U_{ij})(I_{E_j} \ot \wA^{(i)})(u_{i,j}\ot I_\mathcal{K})(I_{E_i}\ot \wA^{(j)^*}).
		\end{equation}
	\end{definition}

   \begin{remark} Equation (\ref{doubly}) admits an equivalent representation given by\begin{align}\label{eqn1.1}
           \wA^{(j)^*} \wA^{(i)} &= (I_{E_j} \ot U_{ij})(I_{E_j} \ot \wA^{(i)})(u_{i,j}\ot I_\mathcal{K})(I_{E_i}\ot \wA^{(j)^*})\\&=(I_{E_j} \ot U_{ij}\wA^{(i)})(u_{i,j}\ot I_\mathcal{K})(I_{E_i} \ot \wA^{(j)^*})\nonumber\\&= (I_{E_j} \ot \wA^{(i)}(I_{E_i} \ot U_{ij}))(u_{i,j}\ot I_\mathcal{K})(I_{E_i} \ot \wA^{(j)^*})\nonumber\\&= (I_{E_j} \ot \wA^{(i)})  (u_{i,j} \ot U_{ij})  (I_{E_i} \ot \wA^{(j)^*}).\nonumber
    \end{align}  
 
       \end{remark}

       	\begin{definition}	
		A closed subspace $\mathcal S$ of $ \mathcal{K}$ is said to be {\it invariant} ({respectively,  \it reducing}) (cf. \cite{SZ08}) for a covariant representation $(\sigma, A^{(1)}, \ldots, A^{(k)})$ on $\mathcal{K},$ if   $ \mathcal{S}$ is $(\sigma, A^{(i)})$-invariant(respectively,   $(\sigma, A^{(i)})$-reducing) subspace for  $i\in J_k.$  Then the natural `restriction' $(\sigma,A^{(1)}, \ldots, A^{(k)})|_\mathcal{S}$ gives a new representation
		of $\be$ on $\mathcal S$, which is known as a {\it summand} of $(\sigma, A^{(1)}, \ldots, A^{(k)}).$	
	\end{definition}

	\begin{notation}
		Let $\beta=\{\beta_1, \dots ,\beta_n \}$ be a non-empty subset of $J_k,$ $\beta^* = J_k\setminus\beta$ which is the compliment set of $\beta$ in $J_k,$
        and let $(\sigma, A^{(1)}, \dots, A^{(k)})$ be a doubly twisted, left-invertible covariant representation of  $\mathbb{E}$ on  $\mathcal{K}.$ Then define
		\begin{align}
			\mathcal{N}_{\beta}:=\bigcap_{i=1}^n (\mathcal{K} \ominus \wt{A}^{(\beta_i)}(E_i \otimes \mathcal{K})).
		\end{align}
		If $\beta=\{i\},$ then we can write $\mathcal{N}_i:=\mathcal{K} \ominus \wt{A}^{(i)}(E_i \otimes \mathcal{K}).$ Thus $$\mathcal{N}_{\beta}=\bigcap_{\beta_i\in\beta}\mathcal{N}_{\beta_i}.$$
        Also, we denote $i^* = \{i\}^* = J_k\setminus \{i\},$ for $i\in J_k.$ 
	\end{notation}
    
The following Lemma is an abstraction of \cite[Lemma 3.5]{Ra22}.
    \begin{lemma}
    $\mathcal{N}_{\beta}$ reduces $U_{ij}$, and $U_{ij}\mathcal{N}_{\beta}=\mathcal{N}_{\beta}$ for each $i\neq j.$
\end{lemma}
   \begin{proof}
    Using $\wA^{(l)} (I_{E_l} \ot U_{ij}) =
		  U_{ij}\wA^{(l)} $, we get $U_{ij}^*\wA^{(l)}=\wA^{(l)} (I_{E_l} \ot U_{ij}^*)$ for each $i\neq j$ and $l.$ Therefore for each $\eta\in\mathcal{N}_{\beta}$ and $l\in\beta$, we get $\wA^{(l)^*}U_{ij}\eta=(I_{E_l} \ot U_{ij})\wA^{(l)^*}\eta=0$, and $\wA^{(l)^*}U_{ij}^*\eta=(I_{E_l} \ot U_{ij}^*)\wA^{(l)^*}\eta=0.$ Therefore $\mathcal{N}_{\beta}$ reduces $U_{ij}.$ 
     \end{proof}
     The following Proposition is an abstraction of \cite[Proposition 2.2]{La22}.
     
     \begin{proposition}\label{Pr}
		Assume that $(\sigma, A^{(1)}, \dots, A^{(k)})$ to be a doubly twisted covariant representation of  $\mathbb{E}$ on  $\mathcal{K}$ such that each $(\sigma, A^{(i)})$ is a near-isometric or satisfies $(1)-(4)$ of Theorem \ref{MT1}, respectively. Then for each non-empty subset $\beta \subseteq J_{k}$, we have the following:
		\begin{enumerate}
			\item   the subspace $\mathcal{N}_{\beta}$ is $(\sigma, A^{(j)})$-reducing, for $j \notin \beta.$
			\item  $ \mathcal{N}_{\beta}\ominus \wt{A}^{(j)}(E_j \otimes \mathcal{N}_{\beta})=\mathcal{N}_{\beta}\cap \mathcal{N}_{j}.$
			\item  $(\sigma ,A^{(j)})$ is a near-isometric or satisfies $(1)-(4)$ of Theorem \ref{MT1}, respectively.
			\end{enumerate}
		\end{proposition}
        
	\begin{proof}
		Let $\emptyset\neq\beta=\{\beta_1, \dots,\beta_n\}\subseteq J_k$ and $j \notin \beta.$ 
		We have to prove that $\mathcal{N}_{\beta}$ is invariant for $(\sigma,{A}^{(j)}),$ let $x\in \mathcal{N}_{\beta}.$ Therefore $x\in \mathcal{N}_{i},$ for each $i\in \beta$ this implies that $\widetilde{A}^{(i)^*}x=0$ for each $i\in \beta.$ Therefore,
		$\wA^{(i)^*} \wA^{(j)}(\eta_j\ot x) =
		(I_{E_i} \ot \wA^{(j)})  (u_{j,i} \ot U_{ji})  (I_{E_j} \ot \wA^{(i)^*})(\eta_j\ot x)=0$ for each $\eta_j\in I_{E_j}.$ It implies $\wA^{(j)}(E_j\ot \mathcal{N}_{\beta})\subseteq$ $\mathcal{N}(\widetilde{A}^{(i)^*})=\mathcal{N}_{i}$ for each $i\in \beta.$ It gives $\wA^{(i)^*} \wA^{(j)}(\eta_j\ot x)=0.$ Therefore, $ \wA^{(j)}(\eta_j\ot x)\in \mathcal{N}_{i}.$ For the other part, 
        $(I_{E_j} \ot \wA^{(i)^*}) \wA^{(j)^*}x =(u_{i,j} \ot U_{ij})
		(I_{E_i} \ot \wA^{(j)^*})\wA^{(i)^*}x=0.$ It deduces $\wA^{(i)^*}x\in$ $\mathcal{N}(I_{E_j} \ot \wA^{(i)^*})=E_j \ot \mathcal{N}(\wA^{(i)^*})=E_j\ot\mathcal{N}_{i}.$ Therefore $\wA^{(j)^*}\mathcal{N}_{i}\subseteq E_j\ot\mathcal{N}_{i}$ for each $i\in\beta.$
        Thus $\mathcal{N}_{\beta}$ is reducing for $(\sigma,{A}^{(j)}).$\\
		To prove (2), let $x\in\mathcal{N}_{\beta}\cap \mathcal{N}_{j}.$ Then $\widetilde{A}^{(j)^*}x.$ It gives $0=\langle \widetilde{A}^{(j)^*}x,h\rangle = \langle x,\widetilde{A}^{(j)}h\rangle$ for each $h\in E_j\ot \mathcal{N}_{\beta}.$ Therefore $x\in \mathcal{N}_{\beta}\ominus \wt{A}^{(j)}(E_j \otimes \mathcal{N}_{\beta}).$ For the other inclusion, let $x\in \mathcal{N}_{\beta}\ominus \wt{A}^{(j)}(E_j \otimes \mathcal{N}_{\beta}),$ then $\langle \widetilde{A}^{(j)^*}x,h\rangle = \langle x,\widetilde{A}^{(j)}h\rangle=0$ for each $h\in E_j\ot \mathcal{N}_{\beta}.$ But $\widetilde{A}^{(j)^*}x\in \mathcal{N}_{\beta},$ as $\mathcal{N}_{\beta}$ is invariant under $\widetilde{A}^{(j)^*}.$ Therefore $\widetilde{A}^{(j)^*}x=0.$ This implies that $x\in\mathcal{N}_{\beta}\cap \mathcal{N}_{j}.$\\
		To prove (3),  using the fact that  $\mathcal{N}_{\beta}$ reduces $(\sigma ,{A}^{(j)})$ and $(\sigma ,{A}^{(j)})$ is a near-isometric or satisfies $(1)-(4)$ of Theorem \ref{MT1}, respectively we get, $(\sigma ,{A}^{(j)})$ is a near-isometric or satisfies $(1)-(4)$ of Theorem \ref{MT1}, respectively; on $\mathcal{N}_{\beta}.$				
	\end{proof}
    For $i\in J_k$
	and $l \in \bn$ define $\wA^{(i)}_l: E_i^{\ot l}\otimes \mathcal{K}\to \mathcal{K}$ by
	\[ \wA^{(i)}_l (\xi_1 \ot \cdots \ot \xi_l \ot h) := A^{(i)} (\xi_1) \cdots A^{(i)}(\xi_l) h,\]
	where $\xi_1, \ldots, \xi_l \in E_{i}, h \in \mathcal{K}$. Then
	\begin{equation}\label{eqnnn}\wA^{(i)}_l=\wA^{(i)}(I_{E_i} \otimes \wA^{(i)}) \cdots (I_{E^{\otimes l-1}_i} \otimes  \wA^{(i)}).\end{equation}

\begin{remark}
 Let $\wA^{(i)}$ be bounded below (or left-invertible) for each $i\in J_k$, and define $L^i:=(\wA^{(i)^*}\wA^{(i)})^{-1}\wA^{(i)^*}$ which is left inverse of $\wA^{(i)}.$ For each $n\in\mathbb{N},$ we further define $L_n^i:\mathcal{K}\to E_i^{\ot n}\ot\mathcal{K}$ by $L_n^i:=(I_{E_i^{\ot n-1}}\ot L^i)(I_{E_i^{\ot n-2}}\ot L^i)\ldots(I_{E_i}\ot L^i)L^i.$ 
\end{remark}

In the following proposition, we shall show that for every doubly twisted left-invertible covariant representation, the left-inverse also satisfies the doubly twisted-type relation. This result will provide an elementary approach to generalize \cite[Remark 3.1]{TV19}.

\begin{proposition}\label{dt 6}
     Assume that $(\sigma, A^{(1)},\ldots,A^{(k)})$ to be a doubly twisted left-invertible covariant representation of $\mathbb{E}$ on $\mathcal{K}.$ Then
    \begin{enumerate}
        \item $(I_{E_j}\ot \wA^{(i)})(u_{i,j} \ot U_{ij}) ((I_{E_i}\ot L^j) = L^j\wA^{(i)}.$
        \item $(I_{E_i}\ot L^j)\wA^{(i)^*} = (u_{j,i} \ot U_{ji})(I_{E_j}\ot \wA^{(i)^*})L^j.$ \label{eqn L_1}
        \item $((I_{E_j}\ot L^i)L^j = (u_{i,j} \ot U_{ij})(I_{E_i}\ot L^j)L^i.$
        \item $(I_{E_i}\ot L^{j^*})(u_{j,i} \ot U_{ji})(I_{E_j}\ot L^i) = L^iL^{j^*}.$
    \end{enumerate}
    \end{proposition}
    \begin{proof}
    $(1)$ By doubly twisted relation (\ref{eqn a}), we obtain
    \begin{align*}
        &  \wA^{(j)^*}\wA^{(i)}(I_{E_i}\ot \wA^{(j)}) = (I_{E_j}\ot \wA^{(i)})(u_{i,j} \ot U_{ij})(I_{E_i}\ot \wA^{(j)^*})(I_{E_i}\ot \wA^{(j)}) .\end{align*}
        It implies \begin{align*}
            \wA^{(j)^*}\wA^{(j)}(I_{E_j}\ot \wA^{(i)})(u_{i,j} \ot U_{ij}) = (I_{E_j}\ot \wA^{(i)})(u_{i,j} \ot U_{ij})((I_{E_i}\ot \wA^{(j)^*}\wA^{(j)}).
        \end{align*} 
         
    Therefore
    \begin{align}\label{li 1}
        & (I_{E_j}\ot \wA^{(i)})(u_{i,j} \ot U_{ij}) ((I_{E_i}\ot (\wA^{(j)^*}\wA^{(j)})^{-1}) = (\wA^{(j)^*}\wA^{(j)})^{-1}(I_{E_j}\ot \wA^{(i)})(u_{i,j} \ot U_{ij}).
    \end{align}It follows
         \begin{align} \label{li 2}
             (I_{E_j}\ot \wA^{(i)})(u_{i,j} \ot U_{ij}) ((I_{E_i}\ot L^j) &= (\wA^{(j)^*}\wA^{(j)})^{-1}(I_{E_j}\ot \wA^{(i)})(u_{i,j} \ot U_{ij})(I_{E_i}\ot \wA^{(j)^*})\nonumber\\& = L^j\wA^{(i)}.
         \end{align}
      $(2)$ By taking the adjoints of both sides of equation (\ref{li 1}), we obtain 
      \begin{align*}
           &  (I_{E_i}\ot (\wA^{(j)^*}\wA^{(j)})^{-1})(u_{j,i} \ot U_{ji})(I_{E_j}\ot \wA^{(i)^*}) = (u_{j,i} \ot U_{ji})(I_{E_j}\ot \wA^{(i)^*})(\wA^{(j)^*}\wA^{(j)})^{-1} \\&
         \implies (I_{E_i}\ot L^j)\wA^{(i)^*} = (u_{j,i} \ot U_{ji})(I_{E_j}\ot \wA^{(i)^*})L^j.
    \end{align*}
     $(3)$ Using the above relation (2), we have 
     \begin{align*}
           &  (u_{j,i} \ot U_{ji})(I_{E_j}\ot \wA^{(i)^*})L^j = (I_{E_i}\ot L^j)\wA^{(i)^*} \\&
         \implies (u_{j,i} \ot U_{ji})(I_{E_j}\ot \wA^{(i)^*})L^j\wA^{(i)} = (I_{E_i}\ot L^j)\wA^{(i)^*}\wA^{(i)} \\&
         \implies (u_{j,i} \ot U_{ji})(I_{E_j}\ot \wA^{(i)^*})(I_{E_j}\ot \wA^{(i)})(u_{i,j} \ot U_{ij}) ((I_{E_i}\ot L^j) = (I_{E_i}\ot L^j)\wA^{(i)^*}\wA^{(i)}.
     \end{align*} Therefore
     \begin{align}\label{li 3}
         (u_{i,j} \ot U_{ij})(I_{E_i}\ot L^j)(\wA^{(i)^*}\wA^{(i)})^{-1} = (I_{E_j}\ot (\wA^{(i)^*}\wA^{(i)})^{-1})(u_{i,j} \ot U_{ij})((I_{E_i}\ot L^j)
     \end{align}
     \begin{align*}
          \implies (u_{i,j} \ot U_{ij})(I_{E_i}\ot L^j)L^i &= (I_{E_j}\ot (\wA^{(i)^*}\wA^{(i)})^{-1})(u_{i,j} \ot U_{ij})((I_{E_i}\ot L^j)\wA^{(i)^*} \\&
        = (I_{E_j}\ot (\wA^{(i)^*}\wA^{(i)})^{-1})(I_{E_j}\ot \wA^{(i)^*})L^j.
     \end{align*} Therefore
     \begin{align}\label{li 4}
         & (u_{i,j} \ot U_{ij})(I_{E_i}\ot L^j)L^i = (I_{E_j}\ot L^i)L^j.
     \end{align}
     $(4)$ By taking the adjoints of both sides of equation (\ref{li 3}) and multiplying corresponding terms, we get 
      \begin{align*}
           &(I_{E_i}\ot L^{j^*})(u_{j,i} \ot U_{ji})(I_{E_j}\ot L^i) = (\wA^{(i)^*}\wA^{(i)})^{-1}(I_{E_i}\ot L^{j^*})(u_{j,i} \ot U_{ji})(I_{E_j}\ot \wA^{(i)^*}) \\& 
              \implies (I_{E_i}\ot L^{j^*})(u_{j,i} \ot U_{ji})(I_{E_j}\ot L^i) = (\wA^{(i)^*}\wA^{(i)})^{-1}(I_{E_i}\ot L^{j^*})(u_{j,i} \ot U_{ji})(I_{E_j}\ot \wA^{(i)^*}) \\& 
               \implies (I_{E_i}\ot L^{j^*})(u_{j,i} \ot U_{ji})(I_{E_j}\ot L^i) = (\wA^{(i)^*}\wA^{(i)})^{-1}\wA^{(i)^*}L^{j^*} = L^iL^{j^*}.\qedhere
      \end{align*}
       \end{proof}

     \begin{remark}\label{RL1}
         By relation (\ref{eqn L_1}) in Proposition \ref{dt 6}, we have $(I_{E_j}\ot \wA^{(i)^*})L^j = (u_{i,j} \ot U_{ij})(I_{E_i}\ot L^j)\wA^{(i)^*}.$ Then 
         \begin{align*}
             &(I_{E_j^{\ot 2}}\ot \wA^{(i)^*})L_2^j \\& = (I_{E_j^{\ot 2}}\ot \wA^{(i)^*})(I_{E_j}\ot L^j)L^j
             = [I_{E_j}\ot(I_{E_j}\ot \wA^{(i)^*})L^j]L^j
             = [I_{E_j}\ot(u_{i,j} \ot U_{ij})(I_{E_i}\ot L^j)\wA^{(i)^*}]L^j
             \\&= (I_{E_j}\ot u_{i,j} \ot U_{ij})(I_{E_j\ot E_i}\ot L^j)(I_{E_j}\ot \wA^{(i)^*})L^j\\&
             = (I_{E_j}\ot u_{i,j} \ot U_{ij})(I_{E_j\ot E_i}\ot L^j)(u_{i,j} \ot U_{ij})(I_{E_i}\ot L^j)\wA^{(i)^*}\\&
             = (I_{E_j}\ot u_{i,j} \ot U_{ij})(u_{i,j} \ot L^j)(I_{E_i\ot E_j}\ot U_{ij})(I_{E_i}\ot L^j)\wA^{(i)^*}\\&
             = (I_{E_j}\ot u_{i,j} \ot U_{ij})( u_{i,j}\ot I_{E_j} \ot U_{ij})(I_{E_i\ot E_j}\ot L^j)(I_{E_i}\ot L^j)\wA^{(i)^*}\\&
             = (I_{E_j}\ot u_{i,j} \ot U_{ij})( u_{i,j}\ot I_{E_j} \ot U_{ij})(I_{E_i}\ot L_2^j)\wA^{(i)^*}.
         \end{align*}
         Similarly, for each $l\in \mathbb{N},$ we have $$(I_{E_j^{\ot l}} \ot\wA^{(i)^*})L_l^j = [\Pi_{r=1}^l(I_{E_j^{\ot l-r}}\ot u_{ij}\ot I_{E_j^{\ot r}}\ot U_{ij})](I_{E_j}\ot L_l^j)\wA^{(i)^*}.$$
         
     \end{remark}
      
      The following Corollary is an abstraction of \cite[Corollary 3.2]{NZ25} and \cite[Theorem 4.9]{HS19} and the proof is similar.
      
     \begin{corollary}
     Assume that $(\sigma, A^{(1)},\ldots,A^{(k)})$ is a left-invertible twisted covariant representation of  $\mathbb{E}$ on  $\mathcal{K}$ such that each $(\sigma, A^{(i)})$ with $i\in J_k$ satisfies any one of the conditions $(1)-(4)$ of Theorem \ref{MT1} or a near-isometric covariant representation. Then 
    \begin{enumerate}
        \item the representation $(\sigma, A^{(1)},\ldots,A^{(k)})$ is doubly twisted, and 
        \item $\wA^{(i)}$ is pure for every $i\in J_k$
    \end{enumerate}
    if and only if 
    \begin{itemize}
        \item[(a)] $(\sigma, A^{(1)},\ldots,A^{(k)})$ has the generating wandering subspace property,
        \item[(b)] $\wA^{(j)}$ and $\wA^{(i)^*}\wA^{(i)}$ satisfies the twisted commutation relation for each $i<j,$ and
        \item[(c)] $\bigvee_{\textbf{m}\in \mathbb{Z}_+^{k-1}} \mathfrak{A}_\textbf{m}^{(i^*)}(\mathcal{N}_{i^*}) \subset \mathcal{N}_i.$
    \end{itemize}
\end{corollary}
The following Lemma is based on \cite[Lemma 1.3]{SZ08}.

\begin{lemma}
    Suppose that $(\sigma, A^{(1)},\ldots,A^{(k)})$ is a twisted left-invertible covariant representation of $\mathbb{E}$ on $\mathcal{K}.$ Then for each $i\in J_k$ the  map
         $\chi^i:\sigma(\mathcal{A})'\to \sigma(\mathcal{A})'$ defined by $\chi^i(x):=\wA^{(i)}(I_{E_i} \otimes x)L^i$ is a normal endomorphism for each $x\in \sigma(\mathcal{A})'.$ Moreover, $\chi^i_l(x)=\wA_l^{(i)}(I_{E_i^{\ot l}} \otimes x)L_l^i$ for $l\in \mathbb{N}$ and for each $x\in \sigma(\mathcal{A})'$, $\chi^i(\chi^j(x)) = \chi^j(\chi^i(x)).$ 
    \end{lemma}

 \begin{proof} For all $x,y\in \sigma(\mathcal{A})',$ we have 
    \begin{align*}
         \chi^i(xy) &= \wA^{(i)}(I_{E_i} \otimes xy)L^i = \wA^{(i)}(I_{E_i} \otimes x)(I_{E_i} \otimes y)L^i \\&= \wA^{(i)}(I_{E_i} \otimes x)L^i\wA^{(i)}(I_{E_i} \otimes y)L^i = \chi^i(x)\chi^i(y).
    \end{align*}
      Now using equations (\ref{li 4}) and (\ref{eqn a}), we get
      \begin{align*}
         \chi^i(\chi^j(x)) & = \wA^{(i)}(I_{E_i} \otimes \chi^j(x))L^i  =  \wA^{(i)}(I_{E_i} \otimes \wA^{(j)}(I_{E_j} \otimes x)L^j)L^i\\ & = \wA^{(j)}(I_{E_j} \otimes \wA^{(i)})(u_{i,j}\ot U_{ij})(I_{E_i\ot E_j}\ot x) (u_{j,i}\ot U_{ji})(I_{E_i} \otimes L^j)L^i \\ &  
        = \wA^{(j)}(I_{E_j} \otimes \wA^{(i)})(I_{E_j\ot E_i}\ot x)(I_{E_i} \otimes L^j)L^i \\ &
        = \wA^{(j)}(I_{E_j} \otimes \wA^{(i)}(I_{E_i} \otimes x)L^i)L^j 
        =  \wA^{(j)}(I_{E_j} \otimes \chi^i(x))L^j 
        = \chi^j(\chi^i(x)). \qedhere
    \end{align*}
    \end{proof}

\begin{remark}
     Using the previous lemma, we acquire the action of $\mathbb{Z}_+^k$ on $\sigma(\mathcal{A})'$ by normal endomorphisms: $\chi(\textbf{n}):= \chi^1_{n_1}\ldots  \chi^k_{n_k}$ (with $\chi(0) = I_{\sigma(\mathcal{A})'})$ for all $\textbf{n}= (n_1, \ldots , n_k)\in \mathbb{Z}_+^k.$ Now we define  $P_l^i: = \chi_l^i(I_\mathcal{K})$ for $l\in \mathbb{N}.$ Further, we define $P(\textbf{n}) = \chi(\textbf{n})(I_{\mathcal{K}})$ (with $P(0)= I_{\mathcal{K}}).$ Denote $P_1^i:=\wA^{(i)}L^i$ and $Q_1^i:=I_\mathcal{K}-P_1^i.$ Then $Q_1^i$ and $P_1^i$ are orthogonal projections and $\text{ran}( Q_1^i)=ker P_1^i = \mathcal{N}_i .$ Therefore, $P_1^i$ is the projection onto $\text{ran}(\wA^{(i)})$.
\end{remark}
The following Proposition is an abstraction of \cite[equation 2.3]{SZ08} for the case of doubly twisted left-invertible covariant representations.

\begin{proposition}
    $P(\textbf{m})P(\textbf{n}) = P(\textbf{m}\vee \textbf{n})$ for each $\textbf{m,n} \in \mathbb{Z}_+^k.$
    \end{proposition}
     \begin{proof}   Using equation (\ref{li 2}), we have
         \begin{align*}
          P_1^iP_1^j &=\wA^{(i)}L^i \wA^{(j)}L^j  = \wA^{(i)}(I_{E_i} \otimes \wA^{(j)})(u_{j,i}\ot U_{ji})(I_{E_j} \otimes L^i)L^j \\& 
          = \wA^{(i)}(I_{E_i} \otimes \wA^{(j)})(I_{E_i} \otimes L^j)L^i 
          = \wA^{(i)}(I_{E_i} \otimes \wA^{(j)}L^j)L^i\\& = \wA^{(i)}(I_{E_i} \otimes P_1^j)L^i  = \chi^i(P_1^j).
      \end{align*} 
      Further, using equations (\ref{li 3}) and (\ref{eqn a}), we get 
      \begin{align*}
           P_1^iP_{l+1}^j=\wA^{(i)}L^i\chi^j(P_l^j) &
          = \wA^{(i)}L^i\wA^{(j)}(I_{E_j} \otimes P_l^j)L^j \\& 
          = \wA^{(i)}(I_{E_i}\otimes \wA^{(j)})(u_{j,i}\ot U_{j,i})(I_{E_j} \otimes L^i)(I_{E_j} \otimes P_l^j)L^j \\&
          = \wA^{(j)}(I_{E_j}\otimes \wA^{(i)})(I_{E_j} \otimes L^iP_l^j)L^j 
          = \wA^{(j)}(I_{E_j} \otimes P_1^iP_l^j)L^j \\&
          = \wA^{(j)}(I_{E_j} \otimes P(e_i+le_j))L^j = \chi^j(P(e_i+le_j)) = P(e_i+(l+1)e_j).
     \end{align*} 
     We can easily check $P_{m_i}^iP_{m_j}^j=P(m_ie_i+m_je_j)$ and therefore $P_{m_j}^jP_{m_i}^i=P_{m_i}^iP_{m_j}^j$ for $m_i,m_j\in \mathbb{N}.$ Further for $i\in J_k, {m_j} , {m_i}\in \mathbb{N}, {m_j}\le {m_i},$ we have $$P_{m_i}^iP_{m_j}^i =  \chi^i_{m_i}(I_{\mathcal{K}})\chi^i_{m_j}(I_{\mathcal{K}}) = \chi^i_{m_j}(\chi^i_{{m_i}-{m_j}}(I_{\mathcal{K}}))\chi^i_{m_j}(I_{\mathcal{K}}) = \chi^i_{m_j}(\chi^i_{{m_i}-{m_j}}(I_{\mathcal{K}})) = \chi^i_{m_i}(I_{\mathcal{K}}) = P_{{m_j}\vee {m_i}}^i.$$ 
    \end{proof}
The following Remark is inspired from \cite[Remark 3.1]{TV19}.

\begin{remark}\label{PC1}
   Let $P^i$ denotes the projection onto $\bigcap_{m_i=1}^\infty\text{ran}(\wA_{m_i}^{(i)})$ which is the SOT-limit of the decreasing sequence of the projections $\{P_{m_i}^i:m_i\in \mathbb{N}\}.$ 
    For all $1\leq i\leq k,$ if $(\sigma,\wA^{(i)})$ admits the Wold-type decomposition(that is, if each $(\sigma,\wA^{(i)})$ satisfies hypothesis of Theorem \ref{MT1} or Theorem \ref{ThmW}), then $Q^i=I-P^i$ is the projection onto $\sum_{m_i=0}^\infty  \mathfrak{A}_{m_i}^{(i)}(\mathcal{N}_i)$ and we get $P^iQ^j=Q^jP^i$ and hence
\begin{align}\label{obs1}
\bigcap_{m_j=0}^\infty\wA_{m_j}^{(j)}(E_j^{\ot m_j}\ot(\sum_{m_i=0}^\infty  \wA_{m_i}^{(i)}(E_i^{\ot m_i}\ot \mathcal{N}_i)))=\sum_{m_i=0}^\infty \wA_{m_i}^{(i)}(E_i^{\ot m_i}\ot(\bigcap_{m_j=0}^\infty \wA_{m_j}^{(j)}(E_j^{\ot m_j}\ot \mathcal{N}_i))),
\end{align}
where $\sum$ is used for $\bigoplus$ or $\bigvee$, respectively.
\end{remark}

\section{Wold-type decomposition for doubly twisted left-invertible covariant representations of product systems}\label{secwd}

In this section, we give a generalized version of \cite[Theorem 5.1]{NZ25} for doubly commuting left-invertible covariant operators to the extended setting of doubly twisted left-invertible covariant representations having Wold-type decomposition. Our approach here is based on \cite{VR17, NZ25, SZ08, Ra22, TV19, HS19}.

\begin{definition}(see the Introduction part of \cite{NZ25}).
    Let $(\sigma,A)$ be a completely bounded covariant representation  of $E$ on $\mathcal{K}.$ If the restriction to any closed invariant subspace has the generating wandering subspace property, then the representation $(\sigma,A)$ is said to admit the {\it Beurling-type theorem}. Equivalently, $\mathcal{S} = \bigvee \mathfrak{A}_n( \mathcal{S}\ominus \wA(E\ot \mathcal{S}))$ for every closed $(\sigma, A)$-invariant subspace $\mathcal{S}.$ 
\end{definition}

The following three Lemmas generalize \cite[Lemmas 2.1-2.3]{NZ25}.
\begin{lemma}\label{BTL1}
    Suppose that $(\sigma, A)$ is a completely bounded covariant representation of $E$ on $\mathcal{K},$ and $\mathcal{S}$ be a reducing subspace for $(\sigma, A).$  Then $(\sigma, A)|_{\mathcal{S}}$ has the generating wandering subspace property and $\mathcal{S} = \bigvee \mathfrak{A}_n( \mathcal{N}(\mathcal{S}^*)),$ if $(\sigma, A)$ has the generating wandering subspace property.
    \end{lemma}
    \begin{proof}
        For a $(\sigma, A)$-reducing subspace $\mathcal{S},$ we get $\wA = \begin{bmatrix}
            \wA|_{E\ot \mathcal{S}} & 0 \\
            0 & \wA|_{E\ot \mathcal{S}^\perp} \\
        \end{bmatrix}.$ We deduce that $(\wA|_{E\ot \mathcal{S}})^* = \wA^*|_{ \mathcal{S}}$ and $\mathcal{N}(\wA|_{E\ot \mathcal{S}})^* = \mathcal{N}(\wA^*|_{ \mathcal{S}}) = \mathcal{N}(\wA^*)\cap \mathcal{S}.$ From the expressions $\mathcal{N}(\wA^*) = \mathcal{N}(\wA^*|_{ \mathcal{S}}) \oplus \mathcal{N}(\wA^*|_{ \mathcal{S}^\perp})$ and $\mathcal{N}(\wA^*) = \bigvee \mathfrak{A}_n(\mathcal{N}(\wA^*|_{ \mathcal{S}}) )\oplus \bigvee \mathfrak{A}_n( \mathcal{N}(\wA^*|_{ \mathcal{S}^\perp}) ),$ we obtain that $\wA|_{E\ot \mathcal{S}}$ has the generating wandering subspace property with $\mathcal{N}(\wA^*)\cap\mathcal{S} = \mathcal{K}\ominus \wA(E\ot \mathcal{K})\cap\mathcal{S} = \mathcal{S}\ominus \wA|_{E\ot \mathcal{S}}(E\ot \mathcal{S}) = \mathcal{S}\ominus \wA(E\ot \mathcal{S}).$
     \end{proof}

\begin{lemma}\label{BTL2}
     Suppose that $(\sigma, A)$ is a completely bounded covariant representation of $E$ on $\mathcal{K},$ and $\mathcal{S}$ be an invariant subspace for $(\sigma, A).$ Then the restriction $(\sigma, A)|_{\mathcal{S}}$ satisfies the Beurling-type theorem, if $(\sigma, A)$ satisfies the Beurling-type theorem.
     \end{lemma}
     \begin{proof}
         Suppose that the $\sigma(\mathcal{A})$-invariant subspace for two closed subspaces $\mathcal{R}$ and $\mathcal{S}$ with $\mathcal{R}\subset \mathcal{S},$ and we denote $\wT = \wA|_{E\ot \mathcal{S}}.$ Therefore, $\wA| _{E\ot \mathcal{R}} = \wT|_{E\ot \mathcal{R}}.$ The remaining part of the proof easily follows by using the Beurling-type theorem for $(\sigma, A).$ 
     \end{proof}

\begin{lemma}\label{wt1}
    Suppose that $(\sigma, A)$ is a completely bounded covariant representation of $E$ on $\mathcal{K}.$ Let $(\sigma, A)$ admits a Wold-type decomposition. Then $(\sigma, A)|_{\mathcal{S}}$ admits a Wold-type decomposition.
\end{lemma}
The following two Lemmas are generalize \cite[Lemma 3.4.2 and 3.4.3]{VR17}.
\begin{lemma}\label{IS1}
    Assume that $(\sigma, A^{(1)},\ldots,A^{(k)})$ to be a doubly twisted left-invertible covariant representation of $\mathbb{E}$ on $\mathcal{K}.$ Then $\mathcal{N}_i$ reduces $L^j$ and $\mathcal{N}_j$ reduces $L^i.$
    \end{lemma}
    \begin{proof}
        By proposition \ref{dt 6}, we have $$L^i\wA^{(j)} = (I_{E_i}\ot \wA^{(j)})(u_{j,i} \ot U_{ji})(I_{E_j}\ot L^i).$$ Then 
        \begin{align*}
            L^i\wA^{(j)}(E_j\ot \mathcal{H}) &= (I_{E_i}\ot \wA^{(j)})(u_{j,i} \ot U_{ji})(I_{E_j}\ot L^i)(E_j\ot \mathcal{H}) \\& 
            = (I_{E_i}\ot \wA^{(j)})(u_{j,i} \ot U_{ji})(E_j\ot L^i\mathcal{H}) 
            \subseteq E_i\ot \wA^{(j)}(E_j\ot \mathcal{H}).
        \end{align*}
        It implies $L^i(\mathcal{N}_j^\perp)\subseteq E_i\ot \mathcal{N}_j^\perp.$
        Again, by Proposition \ref{dt 6}, we have $$(I_{E_j}\ot L^i)\wA^{(j)^*} = (u_{i,j} \ot U_{ij})(I_{E_i}\ot \wA^{(j)^*})L^j.$$ This implies that $\wA^{(j)}(I_{E_j}\ot L^{i^*}) =  L^{i^*}(I_{E_i}\ot \wA^{(j)})(u_{j,i} \ot U_{ji}).$ Thus $L^{i^*}(E_i\ot \mathcal{N}_j^\perp)\subseteq \mathcal{N}_j^\perp.$    
    \end{proof}

\begin{lemma}\label{u1}
    Assume that $(\sigma, A^{(1)},\ldots,A^{(k)})$ to be a doubly twisted left-invertible covariant representation of $\mathbb{E}$ on $\mathcal{K}.$ Then, for each $n\in \mathbb{N},$ $\mathcal{N}_i\cap \text{ran}(\wA_n^{(j)}) = \mathfrak{A}_n^{(j)}( \mathcal{N}_i)$ and $\mathcal{N}_j\cap \text{ran}(\wA_n^{(i)}) = \mathfrak{A}_n^{(i)}( \mathcal{N}_j).$
     \end{lemma}
     \begin{proof}
         We will show the first equation; the second equation is proved in an analogous way. For the reverse inclusion,
         since $\mathcal{N}_i\subseteq \mathcal{K},$ we have $\mathfrak{A}_n^{(j)}( \mathcal{N}_i)\subseteq \text{ran}(\wA_n^{(j)})$ and $\mathcal{N}_i\cap \mathfrak{A}_n^{(j)}( \mathcal{N}_i)\subseteq \mathcal{N}_i\cap \text{ran}(\wA_n^{(j)}).$ Now by lemma \ref{IS1}, $\mathfrak{A}_n^{(j)}( \mathcal{N}_i)\subseteq \mathcal{N}_i$ and therefore, $\mathfrak{A}_n^{(j)}( \mathcal{N}_i)\subseteq \mathcal{N}_i\cap \text{ran}(\wA_n^{(j)}).$ For the direct inclusion, let $x\in \mathcal{N}_i\cap \text{ran}(\wA_n^{(j)}).$ Then, $x\in \mathcal{N}_i$ and $x = \wA_n^{(j)}(\eta_n\ot y),$ for some $\eta_n\ot y\in I_{E_j^{\ot n}}\ot \mathcal{K}.$ We will show that under these conditions $y\in \mathcal{N}_i.$ To see this, for all $z\in \mathcal{N}_i,$ $$\langle y, z\rangle = \langle L_n^{j}\wA_n^{(j)}(\eta_n\ot y), z\rangle = \langle \wA_n^{(j)}(\eta_n\ot y), L_n^{j^*}z\rangle = \langle x, L_n^{j^*}z\rangle.$$ By Lemma \ref{IS1}, $\wA_n^{(j)}(\eta_n\ot z)\in \mathcal{N}_i.$ Therefore, $\langle x, L_n^{j^*}z\rangle = 0$ for all $z\in \mathcal{N}_i.$ Therefore $y \in \mathcal{N}_i,$ as desired. Hence $x\in \mathfrak{A}_n^{(j)}( \mathcal{N}_i).$
     \end{proof}

The following two Lemmas are motivated by \cite[Section 5]{NZ25}.
\begin{lemma}\label{s1}
    Assume that $(\sigma, A^{(1)},\ldots,A^{(k)})$ to be a doubly twisted left-invertible covariant representation of $\mathbb{E}$ on $\mathcal{K}.$ Then $[\bigcap_{n\ge 0} \text{ran}(\wA_n^{(i)})]\cap \mathcal{N}_j = \bigcap_{n\ge 0}\mathfrak{A}_n^{(j)}( \mathcal{N}_i).$ 
     \end{lemma}
     \begin{proof}
      By Lemma \ref{u1}, we have $\mathcal{N}_j\cap \text{ran}(\wA_n^{(i)}) = \mathfrak{A}_n^{(i)}( \mathcal{N}_j).$ Then 
         \begin{align*}
             \bigcap_{n\ge 0}\mathfrak{A}_n^{(i)}( \mathcal{N}_j) = \bigcap_{n\ge 0}[\mathcal{N}_j\cap \text{ran}(\wA_n^{(i)})] = [\bigcap_{n\ge 0} \text{ran}(\wA_n^{(i)})]\cap \mathcal{N}_j.
         \end{align*}
     \end{proof}

\begin{lemma}\label{s2}
    Assume that $(\sigma, A^{(1)},\ldots,A^{(k)})$ to be a doubly twisted left-invertible covariant representation of $\mathbb{E}$ on $\mathcal{K}.$ Then $$\Big[\bigvee_{r \in \mathbb{Z}_+}\mathfrak{A}^{(i)}_{r}(\mathcal{N}_i)\Big] \bigcap \mathcal{N}_j = \bigvee_{r \in \mathbb{Z}_+}\mathfrak{A}^{(i)}_{r}(\mathcal{N}_{\{i,j\}}),$$ where $\mathcal{N}_{\{i,j\}}:= \mathcal{N}_i\cap \mathcal{N}_j.$
    \end{lemma}
    \begin{proof} Observe that 
        $\mathfrak{A}_r^{(i)}( \mathcal{N}_{\{i,j\}})\subseteq \mathfrak{A}_r^{(i)}( \mathcal{N}_j)\subseteq \mathcal{N}_j.$ So, $\bigvee_{r \in \mathbb{Z}_+}\mathfrak{A}_r^{(i)}(\mathcal{N}_{\{i,j\}})\subseteq \mathcal{N}_j$ and $\mathfrak{A}_r^{(i)}( \mathcal{N}_{\{i,j\}})\subseteq \mathfrak{A}_r^{(i)}( \mathcal{N}_i),$ we have $\bigvee_{r \in \mathbb{Z}_+}\mathfrak{A}^{(i)}_{r}(\mathcal{N}_{\{i,j\}})\subseteq \bigvee_{r \in \mathbb{Z}_+}\mathfrak{A}^{(i)}_{r}(\mathcal{N}_{i}).$ Thus, $$\bigvee_{r\in \mathbb{Z}_+}\mathfrak{A}^{(i)}_{r}(\mathcal{N}_{\{i,j\}})\subseteq \Big[\bigvee_{r \in \mathbb{Z}_+}\mathfrak{A}^{(i)}_{r}(\mathcal{N}_i)\Big] \bigcap \mathcal{N}_j.$$ 
        
        For the other inclusion, let $y\in \Big[\bigvee_{r \in \mathbb{Z}_+}\mathfrak{A}^{(i)}_{r}(\mathcal{N}_i)\Big] \bigcap \mathcal{N}_j.$ This implies that $y\in \bigvee_{r \in \mathbb{Z}_+}\mathfrak{A}^{(i)}_{r}(\mathcal{N}_{i})$ and $y\in \mathcal{N}_j.$ Then, there exists a sequence $y_{l}\to y, where $ $y_l = \sum_{t =0}^{N_l}\wA^{(i)}_t(\eta_t\ot y_{t,l}),$ $\eta_t\in E_i^{\ot t},$ and $y_{t,l}\in \mathcal{N}_{i}.$ By using Remark \ref{RL1}, we have
        \begin{align*}
            0 & = [\Pi_{s=1}^{N_l}(I_{E_i^{\ot {N_l}-s}}\ot u_{j,i}\ot I_{E_i^{\ot s}}\ot U_{ji})](I_{E_i}\ot L_{N_l}^i)\wA^{(j)^*}y_l\\&
            = (I_{E_i^{\ot {N_l}}}\ot\wA^{(j)^*})L_{N_l}^iy_l\\&
            = (I_{E_i^{\ot {N_l}}}\ot\wA^{(j)^*})(\eta_{N_l}\ot y_{{N_l},l})\\&
            = (\eta_{N_l}\ot \wA^{(j)^*}y_{{N_l},l}).
        \end{align*}
        This implies that $y_{{N_l},l}\in \mathcal{N}_j.$
        Therefore $y_l\in \bigvee_{r \in \mathbb{Z}_+}\mathfrak{A}_r^{(i)}(\mathcal{N}_{\{i,j\}}).$ Hence, $y\in \bigvee_{r \in \mathbb{Z}_+}\mathfrak{A}^{(i)}_{r}(\mathcal{N}_{\{i,j\}}).$
    \end{proof}

For $\mathbf{m}=(m_1, \cdots, m_k) \in \mathbb{Z}_+^k $, we use notation $\wA_{\mathbf{m}}:\mathbb{E}(\mathbf{m})\otimes \mathcal{K} \longrightarrow \mathcal{K}$ for
$$\wA_{\mathbf{m}}:=\wA^{(1)}_{m_1}\left(I_{E_1^{\otimes m_1}} \otimes\wA^{(2)}_{m_2}\right) \cdots \left(I_{E_1^{\otimes m_1} \otimes \cdots \otimes E_{k-1}^{\otimes {m_{k-1}}}} \otimes\wA^{(k)}_{m_k}\right).$$

	The linear map $A_{\mathbf{m}}: \mathbb{E}(\mathbf{m}) \longrightarrow B(\mathcal{K})$ is defined by $$A_{\mathbf{m}}(\xi)h:=\wA_{\mathbf{m}}(\xi \otimes h), ~ \xi \in \mathbb{E}(\mathbf{m}), h \in \mathcal{K}.$$
We will use  $J_k$ for $\{1,2, \dots ,k\}.$ Let $\beta=\{\beta_1, \dots ,\beta_n \} \subseteq J_k$,   define $$\mathbb{Z}_+^{\beta}:=\{\mathbf{m}=(m_{\beta_1}, \cdots ,m_{\beta_n})\: : \: m_{\beta_j} \in \mathbb{Z}_+, \:1 \leq j \leq n\}.$$ 
For each $\mathbf{m}=(m_{\beta_1}, \cdots ,m_{\beta_n})\in \mathbb{Z}_+^{\beta}$, the map $\wA_{\mathbf{m}}^{\beta}:\mathbb{E}(\mathbf{m})\otimes \mathcal{K} \longrightarrow \mathcal{K}$ is defined by
$$\wA_{\mathbf{m}}^{\beta}:=\wA^{(\beta_1)}_{m_{\beta_1}}\left(I_{E_{\beta_1}^{\otimes m_{\beta_1}}} \otimes\wA^{(\beta_2)}_{m_{\beta_2}}\right) \cdots \left(I_{E_{\beta_1}^{\otimes m_{\beta_1}} \otimes \cdots \otimes E_{\beta_{n-1}}^{\otimes {m_{\beta_{n-1}}}}} \otimes\wA^{(\beta_n)}_{m_{\beta_n}}\right).$$
	
	Moreover, we denote $$\mathfrak{A}_{\mathbf{m}}^{\beta}(\mathcal{S}):=\bigvee \{A_{m_{\beta_1}}^{(\beta_1)}(\eta_{\beta_1}) \cdots A_{m_{\beta_n}}^{(\beta_p)}(\eta_{\beta_n})h\: : \: \:\: \eta_{\beta_j} \in E_{\beta_j}^{\otimes m_{\beta_j}}, 1 \leq j \leq n , h \in \mathcal{S}\},$$ for  a closed   $\sigma(\mathcal{A})$-invariant  subspace  $\mathcal S.$ Clearly, $\mathfrak{A}_{\mathbf{m}}^{\beta}(\mathcal{S})=\overline{\wA_{\mathbf{m}}^{\beta}(\mathbb{E}(\mathbf{m}) \otimes \mathcal{S})}.$ 
	\begin{definition}\begin{enumerate}\item 
			A closed $\sigma(\mathcal{A})$-invariant subspace $\mathcal{S}$ is said to be {\it wandering} for  the covariant representation $(\sigma, A^{(\beta_1)} ,\dots, A^{(\beta_n)})$ if 
			$$\mathcal{S}\perp\mathfrak{A}_{\mathbf{m}}^{\beta}(\mathcal{S}),~\mbox{for each }~\mathbf{m} \in \mathbb{Z}_+^{\beta} \setminus\{0\}.$$
			\item
			The covariant representation $(\sigma, A^{(\beta_1)} ,\dots, A^{(\beta_n)})$  is said to have the {\it "generating wandering subspace  property"} if there exists a wandering subspace $\mathcal{S} \subseteq \mathcal{K}$ for $(\sigma, A^{(\beta_1)} ,\dots, A^{(\beta_n)})$ such that $[\mathcal{K}]_{A_{\beta}}=\mathcal{K}$, that is, $$\mathcal{K}=\bigvee_{\mathbf{m} \in \mathbb{Z}_+^{\beta}}\mathfrak{A}_{\mathbf{m}}^{\beta}(\mathcal{S}).$$
		\end{enumerate}
	\end{definition}
This theorem serves as a generalized version and main results of this article. 
\begin{theorem}\label{MT2}
    If $(\sigma, A^{(1)}, \dots, A^{(k)})$ is a doubly twisted left-invertible covariant representation of $\mathbb{E}$ on $\mathcal{K}$ and for each $i\in J_k,$ if $(\sigma,\wA^{(i)})$ admits the Wold-type decomposition, then for   $2 \leq m \leq k,$ there exists $2^m$  $(\sigma, A^{(1)}, \ldots,	 A^{(m)})$-reducing subspaces $\{\mathcal{K}_{\beta} \: : \: \beta \subseteq J_m \}$ satisfying $\mathcal{K}=\bigoplus_{\beta \subseteq J_m}\mathcal{K}_\beta$   and for each  $\beta=\{i_1, \dots, i_p\} \subseteq J_m: \: (\sigma, A^{(i_1)}, \dots,A^{(i_p)})|_{\mathcal{K}_\beta} $ has  Wold-type decomposition of  $\mathbb{E}_\beta$ over $\mathbb{Z}^\beta_+$ given by the collection of $C^*$-correspondence $\{E_{i_1}, \dots, E_{i_p}\}$  and $(\sigma, A^{(i)})|_{\mathcal{K}_\beta}$ is  induced if $i\in \beta$ and invertible covariant representation when $i \in \beta^*.$ Indeed, the decomposition is unique and
	\begin{align}
		\mathcal{K}_\beta&=\sum_{\mathbf{m} \in\mathbb{Z}_+^{|\beta|}}\mathfrak{A}_{\mathbf{m}}^{\beta}\Big[\bigcap_{\mathbf{i} \in \mathbb{Z}_+^{|\beta^*|}} \mathfrak{A}_{\mathbf{i}}^{\beta^*}(\mathcal{N}_{\beta})\Big] .
	\end{align}
    \end{theorem}
    \begin{proof}
        To prove the claim, we will use induction. We first prove our claim for $k=2.$ Suppose that $(\sigma, A^{(1)}, A^{(2)})$ is doubly twisted left-invertible covariant representation of $\mathbb{E}$ on $\mathcal{K}$ such that $(\sigma, A^{(1)})$ and $(\sigma, A^{(2)})$ admit a Wold-type decomposition. We have 
        \begin{align*}
            \mathcal{K} & = \sum_{m_1 \in \mathbb{Z}_+}\mathfrak{A}^{(1)}_{m_1}(\mathcal{N}_1)\bigoplus \Big[\bigcap_{m_1 \in \mathbb{Z}_+}\mbox{ran}(\widetilde{A}_{m_1}^{(1)})\Big] .
        \end{align*}
        Now by $\wA^{(i)^*}\wA^{(j)} = (I_{E_i}\ot \wA^{(j)})(u_{j,i} \ot U_{ji})(I_{E_j}\ot \wA^{(i)^*}),$ we obtain that $\bigcap_{m_1 \in \mathbb{Z}_+}\mbox{ran}(\widetilde{A}_{m_1}^{(1)})$ and $\sum_{m_1 \in \mathbb{Z}_+}\mathfrak{A}^{(1)}_{m_1}(\mathcal{N}_1)$ are reducing subspaces for $(\sigma, A^{(2)}),$ and so from Lemma \ref{s1} and \ref{s2} we get $$\Big[\bigcap_{m_1 \in \mathbb{Z}_+}\mbox{ran}(\widetilde{A}_{m_1}^{(1)})\Big] \cap \mathcal{N}_2 = \bigcap_{m_1 \in \mathbb{Z}_+}\mathfrak{A}^{(1)}_{m_1}(\mathcal{N}_2)$$ and $\mathcal{N}_2 \bigcap\sum_{m_1 \in \mathbb{Z}_+}\mathfrak{A}^{(1)}_{m_1}(\mathcal{N}_1)  = \sum_{m_1 \in \mathbb{Z}_+}\mathfrak{A}^{(1)}_{m_1}(\mathcal{N}_1\cap \mathcal{N}_2).$  Applying Lemma \ref{wt1} for $(\sigma, A^{(2)})|_{\bigcap_{m_1 \in \mathbb{Z}_+}\mbox{ran}(\widetilde{A}_{m_1}^{(1)})}$ and $(\sigma, A^{(2)})|_{\sum_{m_1 \in \mathbb{Z}_+}\mathfrak{A}^{(1)}_{m_1}(\mathcal{N}_1)}$ and $\widetilde{A}^{(1)}$ is injective, we have
\begin{align*}
             \bigcap_{m_1 \in \mathbb{Z}_+}\mbox{ran}(\widetilde{A}_{m_1}^{(1)}) &= \bigcap_{m_2 \in \mathbb{Z}_+}\mathfrak{A}^{(2)}_{m_2}\Big[\bigcap_{m_1 \in \mathbb{Z}_+}\mbox{ran}(\widetilde{A}_{m_1}^{(1)})\Big]\bigoplus \sum_{m_2 \in \mathbb{Z}_+}\mathfrak{A}^{(2)}_{m_2}\Big[\mathcal{N}_2 \bigcap\big(\bigcap_{m_1 \in \mathbb{Z}_+}\mbox{ran}(\widetilde{A}_{m_1}^{(1)})\big)\Big]\\& 
            = \bigcap_{m_1, m_2 \in \mathbb{Z}_+}\mbox{ran}(\widetilde{A}_{m_1}^{(1)}(I_{E_1^{\ot m_1}}\ot \widetilde{A}_{m_2}^{(2)}))\bigoplus \sum_{m_2 \in \mathbb{Z}_+}\mathfrak{A}^{(2)}_{m_2}\Big[\bigcap_{m_1 \in \mathbb{Z}_+}\mathfrak{A}^{(1)}_{m_1}(\mathcal{N}_2)\Big]
        \end{align*} as well as
        \begin{align*}
             \sum_{m_1 \in \mathbb{Z}_+}\mathfrak{A}^{(1)}_{m_1}(\mathcal{N}_1) &= \bigcap_{m_2 \in \mathbb{Z}_+}\mathfrak{A}^{(2)}_{m_2}\Big[\sum_{m_1 \in \mathbb{Z}_+}\mathfrak{A}^{(1)}_{m_1}(\mathcal{N}_1)\Big]\bigoplus \sum_{m_2 \in \mathbb{Z}_+}\mathfrak{A}^{(2)}_{m_2}\Big[\big(\sum_{m_1 \in \mathbb{Z}_+}\mathfrak{A}^{(1)}_{m_1}(\mathcal{N}_1)\big)\cap \mathcal{N}_2\Big]\\&
            = \bigcap_{m_2 \in \mathbb{Z}_+}\mathfrak{A}^{(2)}_{m_2}\Big[\sum_{m_1 \in \mathbb{Z}_+}\mathfrak{A}^{(1)}_{m_1}(\mathcal{N}_1)\Big]\bigoplus \Big[\sum_{\textbf{m} \in \mathbb{Z}_+^2}\mathfrak{A}^{J_2}_{\textbf{m}}(\mathcal{N}_{J_2})\Big], 
        \end{align*}
        where $\mathcal{N}_{J_2}:= \mathcal{N}_1\cap \mathcal{N}_2$ and $J_2 = \{1,2\}.$
        This implies
    \begin{align*}
    \mathcal{K} = & \mathcal{K}_{\emptyset} \oplus \mathcal{K}_{\{1\}}  \oplus \mathcal{K}_{\{2\}}\oplus \mathcal{K}_{\{1,2\}},
    \end{align*} 
    where 
    \begin{align*}
      & \mathcal{K}_{\emptyset} = \bigcap_{\textbf{m} \in \mathbb{Z}_+^2}\mathfrak{A}^{J_2}_{\textbf{m}}(\mathcal{K}), 
      \mathcal{K}_{\{1\}} = \sum_{m_1 \in \mathbb{Z}_+}\mathfrak{A}^{(1)}_{m_1}\Big[\bigcap_{m_2 \in \mathbb{Z}_+}\mathfrak{A}^{(2)}_{m_2}(\mathcal{N}_1)\Big], \\& 
      \mathcal{K}_{\{2\}} = \sum_{m_2 \in \mathbb{Z}_+}\mathfrak{A}^{(2)}_{m_2}\Big[\bigcap_{m_2 \in \mathbb{Z}_+}\mathfrak{A}^{(1)}_{m_1}(\mathcal{N}_2)\Big], 
      \mathcal{K}_{\{1,2\}} = \sum_{\textbf{m} \in \mathbb{Z}_+^2}\mathfrak{A}^{J_2}_{\textbf{m}}(\mathcal{N}_{J_2}).
  \end{align*}
  Therefore, the statement holds for $n=2.$

  Now, let $l\le n$ and suppose that the statement holds for $(\sigma, A^{(1)}, \dots, A^{(l)}),$ we have $\mathcal{K}=\bigoplus_{\beta \subseteq J_l}\mathcal{K}_\beta,$ where for each nonempty subset of $\beta$ of $J_l,$
  \begin{align}\label{MWD1}
      \mathcal{K}_\beta = \sum_{\mathbf{m} \in\mathbb{Z}_+^p}\mathfrak{A}_{\mathbf{m}}^{\beta}\Big[\bigcap_{\textbf{k} \in \mathbb{Z}_+^{l-p}}\mathfrak{A}^{\beta^*}_{\textbf{k}}(\mathcal{N}_{\beta})\Big] ,
  \end{align} 
  and for $\beta = \emptyset\subseteq J_l,$ $\mathcal{K}_\beta = \bigcap_{\textbf{k} \in \mathbb{Z}_+^{l}}\mathfrak{A}^{\beta}_{\textbf{k}}(\mathcal{K}).$
  By a similar argument as in Remark \ref{PC1}, we have 
   \begin{align}\label{PC2}
       \sum_{\mathbf{m} \in\mathbb{Z}_+^p}\mathfrak{A}_{\mathbf{m}}^{\beta}\Big[\bigcap_{\textbf{k} \in \mathbb{Z}_+^{l-p}}\mathfrak{A}^{\beta^*}_{\textbf{k}}(\mathcal{N}_{\beta})\Big] = \bigcap_{\textbf{k} \in \mathbb{Z}_+^{l-p}}\mathfrak{A}^{\beta^*}_{\textbf{k}}\Big[\sum_{\mathbf{m} \in\mathbb{Z}_+^p}\mathfrak{A}_{\mathbf{m}}^{\beta}(\mathcal{N}_{\beta})\Big].
   \end{align}
   We obtain
   \begin{align*}
       \sum_{\mathbf{m} \in\mathbb{Z}_+^p}\mathfrak{A}_{\mathbf{m}}^{\beta}(\mathcal{N}_{\beta}) = \bigcap_{m_{l+1}\in \mathbb{Z}_+}\mathfrak{A}^{{(l+1)}}_{m_{l+1}}\sum_{\textbf{m} \in \mathbb{Z}_+^p}\mathfrak{A}^{\beta}_{\textbf{m}}(\mathcal{N}_{\beta})\bigoplus \sum_{(r,\textbf{m}) \in \mathbb{Z}_+^{l+1}}\mathfrak{A}^{\beta}_{(r,\textbf{m})}(\mathcal{N}_{\beta\cup \{l+1\}}).
   \end{align*}
   Using equations \ref{MWD1} and \ref{PC2}, we get
   \begin{align*}
       \mathcal{K}_\beta&=\sum_{\mathbf{m} \in\mathbb{Z}_+^p}\mathfrak{A}_{\mathbf{m}}^{\beta}\Big[\bigcap_{\textbf{k} \in \mathbb{Z}_+^{l-p}}\mathfrak{A}^{\beta^*}_{\textbf{k}}(\mathcal{N}_{\beta})\Big]
       = \bigcap_{\textbf{k} \in \mathbb{Z}_+^{l-p}}\mathfrak{A}^{\beta^*}_{\textbf{k}}\Big[\sum_{\mathbf{m} \in\mathbb{Z}_+^p}\mathfrak{A}_{\mathbf{m}}^{\beta}(\mathcal{N}_{\beta})\Big] \\&
       = \bigcap_{\textbf{k} \in \mathbb{Z}_+^{l-p}}\mathfrak{A}^{\beta^*}_{\textbf{k}}\Big[\bigcap_{m_{l+1}\in \mathbb{Z}_+}\mathfrak{A}^{{(l+1)}}_{m_{l+1}}\sum_{\textbf{m} \in \mathbb{Z}_+^p}\mathfrak{A}^{\beta}_{\textbf{m}}(\mathcal{N}_{\beta})\bigoplus \sum_{(r,\textbf{m}) \in \mathbb{Z}_+^{l+1}}\mathfrak{A}^{\beta}_{(r,\textbf{m})}(\mathcal{N}_{\beta\cup \{l+1\}})\Big] \\&
       = \bigcap_{(\textbf{k},r) \in \mathbb{Z}_+^{l-p+1}}\mathfrak{A}^{\beta^*}_{(\textbf{k},r)}\Big[\sum_{\textbf{m} \in \mathbb{Z}_+^p}\mathfrak{A}^{\beta}_{\textbf{m}}(\mathcal{N}_{\beta})\Big] \bigoplus \bigcap_{\textbf{k} \in \mathbb{Z}_+^{l-p}}\mathfrak{A}^{\beta^*}_{\textbf{k}}\Big[\sum_{(r,\textbf{m}) \in \mathbb{Z}_+^{l+1}}\mathfrak{A}^{\beta}_{(r,\textbf{m})}(\mathcal{N}_{\beta\cup \{l+1\}})\Big] \\& 
       = \sum_{\textbf{m} \in \mathbb{Z}_+^p}\mathfrak{A}^{\beta}_{\textbf{m}} \Big[\bigcap_{(\textbf{k},r) \in \mathbb{Z}_+^{l-p+1}}\mathfrak{A}^{\beta^*}_{(\textbf{k},r)}(\mathcal{N}_{\beta})\Big] \bigoplus \sum_{(r,\textbf{m}) \in \mathbb{Z}_+^{l+1}}\mathfrak{A}^{\beta}_{(r,\textbf{m})}\Big[\bigcap_{\textbf{k} \in \mathbb{Z}_+^{l-p}}\mathfrak{A}^{\beta^*}_{\textbf{k}}(\mathcal{N}_{\beta\cup \{l+1\}})\Big].
   \end{align*}
   Applying the Wold-type decomposition to $(\sigma,{A}^{{(l+1)}}),$ we have
   \begin{align*}
       \mathcal{K} & = \sum_{m_{l+1} \in \mathbb{Z}_+}\mathfrak{A}^{{(l+1)}}_{m_{l+1}}(\mathcal{N}_{l+1})\bigoplus \Big[\bigcap_{m_{l+1} \in \mathbb{Z}_+}\mbox{ran}(\widetilde{A}_{{m_{l+1}}}^{{(l+1)}})\Big]
   \end{align*}
   and hence for $\beta = \emptyset \subseteq J_l,$
   \begin{align*}
       \mathcal{K}_\beta &= \bigcap_{\textbf{k} \in \mathbb{Z}_+^{l}}\mathfrak{A}^{\beta}_{\textbf{k}}(\mathcal{K})
       = \bigcap_{\textbf{k} \in \mathbb{Z}_+^{l}}\mathfrak{A}^{\beta}_{\textbf{k}} \Big[\sum_{m_{l+1} \in \mathbb{Z}_+}\mathfrak{A}^{{(l+1)}}_{m_{l+1}}(\mathcal{N}_{l+1})\bigoplus \bigcap_{m_{l+1} \in \mathbb{Z}_+}\mbox{ran}(\widetilde{A}_{{m_{l+1}}}^{{(l+1)}})\Big] \\& 
       = \sum_{m_{l+1} \in \mathbb{Z}_+}\mathfrak{A}^{{(l+1)}}_{m_{l+1}}\Big[\bigcap_{\textbf{k} \in \mathbb{Z}_+^{l}}\mathfrak{A}^{\beta}_{\textbf{k}}(\mathcal{N}_{l+1})\Big]\bigoplus \bigcap_{(\textbf{k},m_{l+1}) \in \mathbb{Z}_+^{l+1}}\mathfrak{A}^{\beta\cup \{l+1\}}_{ (\textbf{k},m_{l+1})}.
   \end{align*}
   Consequently, $\mathcal{K}=\bigoplus_{\beta \subseteq J_{l+1}}\mathcal{K}_\beta.$
    By the above decomposition of $\mathcal{K}$ we conclude that the covariant representation $(\sigma,A^{(r_1)}, \dots, A^{(r_p)})|_{\mathcal{K}_\beta}$ has Wold-type decomposition of $\mathbb{E}_\beta$ over $\mathbb{Z}^\beta_+$  $(\sigma, A^{(i)})|_{\mathcal{K}_\beta}$ is induced for $i\in \beta$ and invertible  for every $i \in  \beta^*.$   
	\end{proof}

The subsequent Theorem is an extended version of \cite[Theorem 5.2]{NZ25}.
\begin{theorem}
    Suppose that $(\sigma, A^{(1)},\ldots,A^{(k)})$ is a doubly twisted left-invertible covariant representation of  $\mathbb{E}$ on $\mathcal{K}$ such that each $(\sigma, A^{(i)})$ admits a Wold-type decomposition. Then we can decompose into two reducing subspaces for $(\sigma, A^{(1)},\ldots,A^{(k)})$ as $\mathcal{K} = \mathcal{K}_1\bigoplus \mathcal{K}_2$ such that $(\sigma, A^{(1)},\ldots,A^{(k)})|_{\mathcal{K}_1}$ is unitary covariant representation and the restriction to $\mathcal{K}_2$ has the generating wandering subspace property. In particular, $(\sigma, A^{(1)},\ldots,A^{(k)})$ admits a von Neumann Wold decomposition. 
     \end{theorem}
    \begin{proof}
        By Theorem \ref{MT2}, we write  $\mathcal{K}=\bigoplus_{\beta \subseteq J_k}\mathcal{K}_\beta.$ For $\mathcal{K}_1 = \mathcal{K}_{\emptyset}$ and $\mathcal{K}_2=\bigoplus_{\emptyset\neq\beta \subseteq J_k}\mathcal{K}_\beta,$ the required results is fulfilled.
    \end{proof}

\begin{notation}(see \cite[section 5] {SS26}).
    For any $n\in \mathbb{Z}_+,$ we can define an isomorphism $u_{i,k}^{(n)}: E_i \ot E_k^{\ot n} \to E_k^{\ot n}\ot E_i$ recursively by 
    \begin{equation}
        u_{i,k}^{(n)}:= \prod_{r=1}^l \big(I_{E_k^{\ot n-r}} \ot u_{i,k}\ot I_{E_k^{\ot r-1}}\big).
    \end{equation}
    Moreover, for each $m,n\in \mathbb{Z}_+,$ we define $u_{i,k}^{(m,n)}: E_i^{\ot m} \ot E_k^{\ot n} \to E_k^{\ot n}\ot E_i^{\ot m}$ recursively by 
    \begin{equation}
        u_{i,k}^{(m,n)}:= \prod_{r=1}^p \big(I_{E_i^{\ot r-1}} \ot u_{i,k}^{(n)}\ot I_{E_k^{\ot p-r}}\big).
    \end{equation}
\end{notation}

The subsequent Lemma is a generalization of \cite[Lemma 4.3]{NS22}.
\begin{lemma}\label{OD1}
     Suppose that $(\sigma, A^{(1)},\ldots,A^{(k)})$ is a twisted left-invertible covariant representation of $\mathbb{E}$ on $\mathcal{K}.$ Then there exists a monomial $\alpha_{i,\textbf{m}}\in \mathbb{C}[\lambda_1,\ldots,\lambda_k]$ such that for each $\textbf{m}\in \mathbb{Z}_+^k,$ and $i\in J_k,$ we get $$\wA^{(i)}(I_{E_i}\ot\wA_\textbf{m}) = \wA_\textbf{m}(I_{\mathbb{E}(\textbf{m})}\ot\wA^{(i)})(u_{i,\textbf{m}}\ot\alpha_{i,\textbf{m}}(U)),$$  where $u_{i,\textbf{m}} = \Pi_{j=1,j\neq i}^k(I_{{E_1}^{\ot m_1}}\ot\ldots \ot I_{E_{j-1}^{\ot m_{j-1}}}\ot u_{i,j}^{(m_j)}\ot  I_{E_{j+1}^{\ot m_{j+1}}}\ot u_{i,j}^{(m_j)}\ldots \ot I_{{E_k}^{\ot m_k}}),$
         $$\alpha_{i,\textbf{m}} = \lambda_1^{m_1}\ldots\lambda_{i-1}^{m_{i-1}}\lambda_{i+1}^{m_{i+1}}\ldots \lambda_{k}^{m_{k}}\in \mathbb{C}[\lambda_1,\ldots,\lambda_k],$$ and $$\alpha_{i,\textbf{m}}(U) = U_{i1}^{m_1}\ldots U_{i(i-1)}^{m_{i-1}}U_{i(i+1)}^{m_{i+1}}\ldots U_{ik}^{m_k}.$$
         \end{lemma}
     \begin{proof}
         Note that $\wA^{(i)} (I_{E_i} \ot \wA^{(j)}) = \wA^{(j)}(I_{E_j}\ot  \wA^{(i)})(u_{i,j}\ot U_{ij}).$ Therefore 
         \begin{align*}
             \wA^{(i)}(I_{E_i}\ot\wA_\textbf{m}) & = \wA_\textbf{m}(I_{\mathbb{E}(\textbf{m})}\ot\wA^{(i)})(u_{i,\textbf{m}}\ot \alpha_{i,\textbf{m}}(U)).\qedhere
         \end{align*} 
         \end{proof}

The Lemma below provides an important abstraction of \cite[Lemma 5.3]{NS22}.
\begin{lemma}\label{OD2}
    Assume that $(\sigma, A^{(1)},\ldots,A^{(k)})$ to be a doubly twisted left-invertible covariant representation of $\mathbb{E}$ on $\mathcal{K},$ and let $\mathcal{S}\subseteq \mathcal{K}$ reduces $\wA_\textbf{m}.$ Then $$\bigcap_{r\in \mathbb{Z}_+}(\wA^{(1)}(I_{E_1}\ot\wA^{(2)})\ldots (I_{E_1^{\ot k-1}}\ot\wA^{(k)}))_r\mathcal{S} = \bigcap_{\textbf{m}\in \mathbb{Z}_+^k}\wA_\textbf{m}({\mathbb{E}(\textbf{m})}\ot\mathcal{S}).$$
    \end{lemma}
    \begin{proof}
        We show here for only $k = 2$ as the remaining part can easily be shown by induction for any positive integer $k\ge 2.$ Suppose that $(\sigma, A^{(1)},A^{(2)})$ is a twisted left-invertible. Define a left-invertible covariant representation $(\sigma, A)$ of $E:=E_1\ot E_2$ on $\mathcal K$ by $$A(\xi_1\ot\xi_2)h:=A^{(1)}(\eta_1)A^{(2)}(\eta_2)h~\mbox{for all}~\xi_1\in E_1,\xi_2\in E_2$$ and $h\in \mathcal{K}.$ Note that $\bigcap_{m_1,m_2\in \mathbb{Z}_+}\wA_{m_1}^{(1)}(I_{E_1^{\ot m_1}}\ot\wA_{m_2}^{(2)})(E_1^{\ot m_1}\ot E_2^{\ot m_2}\ot \mathcal{S}) \subseteq \bigcap_{r\in \mathbb{Z}_+}\wA_r(E^{\ot r}\ot \mathcal{S}).$ By Lemma \ref{OD1}, for all $n \in \mathbb{Z}_+$ we get  
        \begin{align*}
            \wA_n & = \wA_n^{(1)}(I_{E_1^{\ot n}}\ot\wA_n^{(2)})(I_{E_1^{\ot n-1}}\ot u_{1,2}^{(n-1)^*}\ot I_{E_2}\ot U_{12}^{n-1^*})(I_{E_1^{\ot n-2}}\ot u_{1,2}^{(n-2)^*}\ot I_{E_2\ot E}\ot U_{12}^{n-2^*})\ldots \\&\quad \quad \quad (I_{E_1^{\ot 2}}\ot u_{1,2}^{(2)^*}\ot I_{E_2\ot E^{\ot n-3}}\ot U_{12}^{2^*})(I_{E_1}\ot u_{1,2}^{*}\ot I_{E_2\ot E^{\ot n-2}}\ot U_{12}^{*}).
        \end{align*}
         The required other inclusion follows from the above expression.
    \end{proof}

The subsequent Corollary is an abstraction of \cite[Proposition 5.4]{NS22} and it is inspired by weak bishift notion due to Popovici \cite{P04} in the setting of left-invertible covariant representation.
\begin{corollary}\label{OD3}
    Assume that $(\sigma, A^{(1)},\ldots,A^{(k)})$ to be a doubly twisted left-invertible covariant representation of $\mathbb{E}$ on $\mathcal{K}$ such that each $(\sigma, A^{(i)})$ has the Wold-type decomposition. Then $(\sigma, A^{(1)},\ldots,A^{(k)})$ is twisted analytic if and only if $ \wA^{(i)}|_{\mathcal{N}_{i^*}}$ and $\wA^{(j)}(I_{E_j}\ot \wA^{(r)})$ are analytic for $i,j,r\in J_k$ and $j\ne r.$ 
    \end{corollary}
    \begin{proof}
        Suppose that $(\sigma, A^{(1)},\ldots,A^{(k)})$ is an analytic, then as $\mathcal{N}_{i^*}$ reduces $(\sigma, A^{(i)}),$ it follows that $\wA^{(i)}|_{\mathcal{N}_{i^*}}$ is analytic. Moreover, by Lemma \ref{OD2}, we have 
        \begin{align*}
            \bigcap_{m_2\in \mathbb{Z}_+}\text{ran}(\wA_{m_2}) & = \bigcap_{m_j,m_r\in \mathbb{Z}_+}\text{ran}\wA_{m_j}^{(j)}(I_{E_j^{\ot m_j}}\ot\wA_{m_r}^{(r)}) \\& 
            \subseteq \bigcap_{m_j\in \mathbb{Z}_+}\text{ran}\wA_{m_j}^{(j)} = \{0\},
        \end{align*}
        as $(\sigma,A^{(j)})$ is analytic. Therefore, $\wA^{(j)}(I_{E_j}\ot \wA^{(r)})$ is analytic for $j \neq r.$ Conversely, suppose that $ \wA^{(i)}|_{\mathcal{N}_{i^*}}$ and $\wA^{(j)}(I_{E_j}\ot \wA^{(r)})$ are analytic $j\ne r$ and $i\in J_k.$ Our aim is to show that $\wA^{(j)}$ is analytic, that is, $\bigcap_{m_i\ge 0} \text{ran}\wA_{m_i}^{(i)} = 0,$ for $i\in J_k.$ Note that $(\sigma, A^{(1)},\ldots,A^{(i-1)},A^{(i+1)},\ldots,A^{(k)})$ is $(k-1)$ -tuple of doubly twisted left-invertible covariant representation with respect to the twist $\{U_{cd}:c<d,c,d\ne i\}.$ Then by Theorem \ref{MT2}, $(\sigma, A^{(1)},\ldots,A^{(i-1)},A^{(i+1)},\ldots,A^{(k)})$ admits von Neumann-Wold type decomposition $\mathcal{K}=\bigoplus_{\beta \subseteq i^*}\mathcal{K}_\beta,$ where
        \begin{align*}
            \mathcal{K}_\beta&=\sum_{\mathbf{m} \in\mathbb{Z}_+^{|\beta|}}\mathfrak{A}_{\mathbf{m}}^{\beta}\Big[\bigcap_{\mathbf{i} \in \mathbb{Z}_+^{ k-1-|\beta|}} \mathfrak{A}_{\mathbf{i}}^{\beta^*}(\mathcal{N}_{\beta})\Big] \quad \quad \quad \quad (\beta \subseteq i^*).
        \end{align*}
         Fix $\beta\subseteq i^*$ and $m_i\in \mathbb{Z}_+.$ Suppose that $t = k-1-|\beta|.$ By Lemma \ref{OD1}, for each $\textbf{m}\in \mathbb{Z}_+^{|\beta|},$ there exist a product system unitary $u_\textbf{m}$ and a monomial $\alpha_\textbf{m}$ such that $\wA_{m_i}^{(i)}(I_{E_i^{\ot m_i}}\ot\wA_{\textbf{m}}^{(\beta)}) = \wA_{\textbf{m}}^{(\beta)}(I_{E_\beta^{\ot \textbf{m}}}\ot\wA_{m_i}^{(i)})(u_\textbf{m}\ot \alpha_\textbf{m}(U))$. We obtain
        \begin{align*}
            \wA_{m_i}^{(i)}(I_{E_i^{\ot m_i}}\ot\mathcal{K}_\beta) & = \mathfrak{A}_{m_i}^{(i)}\left(\sum_{\textbf{m} \in\mathbb{Z}_+^{|\beta|}}\mathfrak{A}_{\mathbf{m}}^{\beta}\Big[\bigcap_{\mathbf{i} \in \mathbb{Z}_+^{t}} \mathfrak{A}_{\textbf{i}}^{i^*\setminus\beta}(\mathcal{N}_{\beta})\Big]\right) \\&
            = \sum_{\textbf{m} \in\mathbb{Z}_+^{|\beta|}}\wA_{\textbf{m}}^{(\beta)}(I_{E_\beta^{\ot \textbf{m}}}\ot\wA_{m_i}^{(i)})(u_\textbf{m}\ot \alpha_\textbf{m}(U)))\Big[\bigcap_{\textbf{i} \in \mathbb{Z}_+^{t}} \mathfrak{A}_{\mathbf{i}}^{i^*\setminus\beta}(\mathcal{N}_{\beta})\Big] \\&
            = \sum_{\textbf{m} \in\mathbb{Z}_+^{|\beta|}}\wA_{\textbf{m}}^{(\beta)}(I_{E_\beta^{\ot \textbf{m}}}\ot\wA_{m_i}^{(i)})\Big[\bigcap_{\textbf{i} \in \mathbb{Z}_+^{t}} \mathfrak{A}_{\mathbf{i}}^{i^*\setminus\beta}(\mathcal{N}_{\beta})\Big] \\&
            = \sum_{\textbf{m} \in\mathbb{Z}_+^{|\beta|}}\wA_{\textbf{m}}^{(\beta)}\Big[\bigcap_{\textbf{i} \in \mathbb{Z}_+^{t}}(I_{E_\beta^{\ot \textbf{m}}}\ot\wA_{m_i}^{(i)}) \mathfrak{A}_{\textbf{i}}^{i^*\setminus\beta}(\mathcal{N}_{\beta})\Big],
        \end{align*}
        where the last but one equality follows from $$(u_\textbf{m}\ot \alpha_\textbf{m}(U))\Big[\bigcap_{\mathbf{i} \in \mathbb{Z}_+^{t}} \mathfrak{A}_{\textbf{i}}^{i^*\setminus\beta}(\mathcal{N}_{\beta})\Big] = \Big[\bigcap_{\textbf{i} \in \mathbb{Z}_+^{t}} \mathfrak{A}_{\textbf{i}}^{i^*\setminus\beta}(\mathcal{N}_{\beta})\Big].$$
        Therefore $$\bigcap_{m_i\in\mathbb{Z}_+} \wA_{m_i}^{(i)}(I_{E_i^{\ot m_i}}\ot\mathcal{K}_\beta) = \sum_{\textbf{m} \in\mathbb{Z}_+^{|\beta|}}\wA_{\textbf{m}}^{(\beta)}\Big[\bigcap_{\textbf{i} \in \mathbb{Z}_+^{t}}(I_{E_\beta^{\ot \textbf{m}}}\ot\wA_{m_i}^{(i)}) \mathfrak{A}_{\textbf{i}}^{i^*\setminus\beta}(\mathcal{N}_{\beta})\Big].$$
        If $\beta\subset i^*,$ there exists $j\in i^*\setminus\beta,$ such that 
        \begin{equation}\label{ODE1}
            \bigcap_{\textbf{i} \in \mathbb{Z}_+^{t},m_i\in \mathbb{Z}_+}\mathfrak{A}_{m_i}^{(\textbf{i})}(I_{E_i^{\ot m_i}}\ot \mathfrak{A}_{\textbf{i}}^{i^*\setminus\beta}(\mathcal{N}_{\beta})) = \bigcap_{\mathbf{i}' \in \mathbb{Z}_+^{t-1},m_i,m_j\in \mathbb{Z}_+}\wA_{m_i}^{(i)}(I_{E_i^{\ot m_i}}\ot \wA_{m_j}^{(j)})(I_{E_i^{\ot m_i}\ot E_j^{\ot m_j}}\ot \mathfrak{A}_{\mathbf{i}'}^{i^*\setminus\beta\cup \{j\}}(\mathcal{N}_{\beta})),
        \end{equation}
        where $\textbf{m} = (m_1,\ldots,m_k)\in \mathbb{Z}_+^k.$ Applying Lemma \ref{OD2}, we have
        \begin{align*}
            \bigcap_{\mathbf{i}' \in \mathbb{Z}_+^{t},m_i,m_j\in \mathbb{Z}_+}\wA_{m_i}^{(i)}(I_{E_i^{\ot m_i}}\ot \wA_{m_j}^{(j)})(I_{E_i^{\ot m_i}\ot E_j^{\ot m_j}}\ot \mathfrak{A}_{\textbf{i}'}^{i^*\setminus\beta\cup \{j\}}(\mathcal{N}_{\beta})) &\subseteq \bigcap_{m_i,m_j\in \mathbb{Z}_+}\text{ran}\wA_{m_i}^{(i)}(I_{E_i^{\ot m_i}}\ot\wA_{m_j}^{(j)}) \\& 
            \subseteq \bigcap_{m_2\in \mathbb{Z}_+}\text{ran}(\wA_{m_2}) = \{0\}.
        \end{align*}
        Then equation (\ref{ODE1}) implies $\bigcap_{\textbf{i} \in \mathbb{Z}_+^{t-1},m_i\in \mathbb{Z}_+}\wA_{m_i}^{(i)}(I_{E_i^{\ot m_i}}\ot \mathfrak{A}_{\mathbf{i}}^{i^*\setminus\beta}(\mathcal{N}_{\beta})) = \{0\}.$ Then the equality preceding equation (\ref{ODE1}) yields $\bigcap_{_i\in \mathbb{Z}_+}\wA_{m_i}^{(i)}(I_{E_i^{\ot m_i}}\ot \mathcal{K}_\beta) = \{0\}$ for all $\beta\subset i^*.$ Therefore 
        \begin{align*}
            \bigcap_{m_i\in \mathbb{Z}_+}\text{ran}\wA_{m_i}^{(i)}& = \bigcap_{m_i\in \mathbb{Z}_+}\wA_{m_i}^{(i)}(I_{E_i^{\ot m_i}}\ot \mathcal{K}_{i^*})
            = \sum_{\textbf{i}'\in \mathbb{Z}_+^{k-1}}\mathfrak{A}_{\mathbf{i}'}^{i^*}\big(\bigcap_{m_i\in \mathbb{Z}_+}\mathfrak{A}_{m_i}^{i}(\mathcal{N}_{i^*})\big).
        \end{align*}
        As $\wA^{(i)}|_{\mathcal{N}_{i^*}}$ is analytic by assumption, we obtain $\bigcap_{m_i\in \mathbb{Z}_+}\text{ran}\wA_{m_i}^{(i)} = 0.$
    \end{proof}

\paragraph{$\mathbf{Acknowledgment}$}: Niraj Kumar is supported by a UGC fellowship (NTA Ref. No.:231610038861). Harsh Trivedi thanks Science and Engineering Research Board(SERB), Department of Science \& Technology (DST), Government of India for the research grant received through MATRICS-SERB, File No: MTR/2021/000286. We also acknowledge the Centre for Mathematical \& Financial Computing and the DST-FIST grant for financial support for the computing lab facility through the scheme FIST(File No: SR/FST/MS-I/2018/24) at the LNM institute of information technology, Jaipur.

\end{document}